\documentclass[a4paper,11pt]{article}
\usepackage[top=3cm,bottom=3cm,left=2.5cm,right=2.5cm,a4paper]{geometry}
\usepackage{amsfonts,amsmath,amssymb,amsthm}
\usepackage{float,xcolor,fancyhdr,siunitx}
\usepackage{setspace}
\usepackage{lineno}
\usepackage[normalem]{ulem} 
\usepackage[colorlinks]{hyperref}
\hypersetup{
    colorlinks=true, 
    linkcolor={blue!50!black},
    citecolor={green!60!black},
    urlcolor= {red!80!black}
}
\usepackage{enumitem,mdframed}
\usepackage[]{algorithm2e}
\usepackage{empheq}
\usepackage[bottom]{footmisc} 
\usepackage{tikz,pgfplots}
\usepgfplotslibrary{groupplots}
\usetikzlibrary{positioning,arrows,fit,patterns}
\pgfplotsset{compat=newest}
\usepackage[labelformat=simple]{subcaption}       
\DeclareCaptionLabelFormat{opening}{(#2)} 
\makeatletter
\renewcommand*\env@matrix[1][\arraystretch]{%
  \edef\arraystretch{#1}%
  \hskip -\arraycolsep
  \let\@ifnextchar\new@ifnextchar
  \array{*\c@MaxMatrixCols c}}
\makeatother

\newcommand{\update}[1]{#1}

\newcommand*\samethanks[1][\value{footnote}]{\footnotemark[#1]}

\newcommand{\forward}     {\mathcal{F}}
\newcommand{\R}           {\mathbb{R}}
\newcommand{\Real}        {\mathrm{Re}}
\newcommand{\Cx}          {\mathbb{C}}

\newcommand{\divergence}  {\nabla \cdot}
\newcommand{\ii}          {\mathrm{i}}
\newcommand{\dd}          {\mathrm{d}}
\newcommand{\bx}          {\boldsymbol{x}}
\newcommand{\misfit}      {\mathcal{J}}
\newcommand{\pressure}    {p}
\newcommand{\velocity}    {\boldsymbol{v}}
\newcommand{\freq}        {\sigma}

\newcommand{\srcp}        {f}
\newcommand{\n}           {\boldsymbol{\nu}} 
\newcommand{\vn}          {\velocity_{\n}}   
\newcommand{\data}        {\boldsymbol{d}} 
\newcommand{\model}       {\boldsymbol{m}}

\newcommand{\psitest}     {{\boldsymbol{\psi}}}
\newcommand{\psitestc}    {\overline{\boldsymbol{\psi}}}

\newcommand{\phic}        {\overline{\phi}}
\newcommand{\xic}         {\overline{\xi}}
\newcommand{\mesh}      {{\mathcal{T}_h}}
\newcommand{\facesetin} {\Sigma^{\text{I}}}
\newcommand{\facesetout}{\Sigma^{\text{B}}}
\newcommand{\facesetall}{\Sigma}           
\newcommand{\nfaceall}  {{N_{\facesetall}}} 
\newcommand{\nfacein}   {{N_{\facesetin}}}  
\newcommand{\nfaceout}  {{N_{\facesetout}}} 

\newcommand{\face}     {\mathfrak{f}}
\newcommand{\ncell}    {{N}}
\newcommand{\cell}     {{K_e}}
\newcommand{\dcell}    {{\partial\cell}}
\newcommand{\hdgA}{\mathbb{A}}
\newcommand{\hdgB}{\mathbb{B}}
\newcommand{\hdgC}{\mathbb{C}}
\newcommand{\hdgL}{\mathbb{L}}
\newcommand{\hdgS}{\mathbb{S}}

\newcommand{\hdgProject}{\mathcal{R}}
\newcommand{\lagrangian}{\mathcal{L}}
\newcommand{\tildu}     {{\tilde{U}}}
\newcommand{\tildL}     {\tilde{\Lambda}}
\newcommand{\tildadjA}  {\tilde{\gamma}_{1,e}}
\newcommand{\tildadjB}  {\tilde{\gamma}_2}
\newcommand{\adjA}  {\gamma_{1,e}}
\newcommand{\adjB}  {\gamma_2}
\newlength{\modelwidth} \newlength{\modelheight}
\newcommand{\modelfile} {}
 
\newlength {\stepalgo}

\newtheorem{remark}     {Remark}

\usepackage[capitalize,nameinlink]{cleveref}
\crefname{section}   {Section}   {Sections}
\crefname{subsection}{Subsection}{Subsections}
\Crefname{section}   {Section}   {Sections}
\Crefname{subsection}{Subsection}{Subsections}
\Crefname{figure}    {Figure}    {Figures}

\crefformat{equation}{\textup{#2(#1)#3}}
\crefrangeformat{equation}{\textup{#3(#1)#4--#5(#2)#6}}
\crefmultiformat{equation}{\textup{#2(#1)#3}}{ and \textup{#2(#1)#3}}
{, \textup{#2(#1)#3}}{, and \textup{#2(#1)#3}}
\crefrangemultiformat{equation}{\textup{#3(#1)#4--#5(#2)#6}}%
{ and \textup{#3(#1)#4--#5(#2)#6}}{, \textup{#3(#1)#4--#5(#2)#6}}{, and \textup{#3(#1)#4--#5(#2)#6}}
\Crefformat{equation}{#2Equation~\textup{(#1)}#3}
\Crefrangeformat{equation}{Equations~\textup{#3(#1)#4--#5(#2)#6}}
\Crefmultiformat{equation}{Equations~\textup{#2(#1)#3}}{ and \textup{#2(#1)#3}}
{, \textup{#2(#1)#3}}{, and \textup{#2(#1)#3}}
\Crefrangemultiformat{equation}{Equations~\textup{#3(#1)#4--#5(#2)#6}}%
{ and \textup{#3(#1)#4--#5(#2)#6}}{, \textup{#3(#1)#4--#5(#2)#6}}{, and \textup{#3(#1)#4--#5(#2)#6}}

\crefdefaultlabelformat{#2\textup{#1}#3}

\crefname{proposition}{Proposition}{Propositions}
\Crefname{proposition}{Proposition}{Propositions}
\crefname{definition} {Definition} {Definitions}
\Crefname{definition} {Definition} {Definitions}
\crefname{theorem}    {Theorem}    {Theorems}
\Crefname{theorem}    {Theorem}    {Theorems}
\crefname{remark}     {Remark}     {Remarks}
\Crefname{remark}     {Remark}     {Remarks}
\crefname{assumption} {Assumption} {Assumptions}
\Crefname{assumption} {Assumption} {Assumptions}

\title{Adjoint-state method for Hybridizable Discontinuous Galerkin discretization: 
       application to the inverse acoustic wave problem}

\author{
Florian Faucher\thanks{Faculty of Mathematics, University of Vienna, Oskar-Morgenstern-Platz 1,
                       A-1090 Vienna, Austria.
                      (\href{mailto:florian.faucher@univie.ac.at}
                      {\texttt{florian.faucher@univie.ac.at}}).} 
\and
Otmar Scherzer\samethanks[1]
              \thanks{Johann Radon Institute for Computational and 
                      Applied Mathematics (RICAM), Altenbergerstra{\ss}e 69 A-4040, Linz, Austria.}
}
\date{}

\begin{document}
\maketitle 

\begin{abstract}
  In this paper, we perform non-linear minimization using 
  the Hybridizable Discontinuous Galerkin method (HDG) for
  the discretization of the forward problem, and implement 
  the adjoint-state method for the efficient computation of 
  the functional derivatives.
  Compared to continuous and discontinuous Galerkin discretizations, 
  HDG reduces the computational cost by using the
  numerical traces for the global linear system, 
  hence removing the degrees of freedom that are 
  inside the cells.
  It is particularly attractive for large-scale time-harmonic 
  quantitative inverse problems which make repeated use of the 
  forward discretization as they rely on an iterative minimization 
  procedure.
  HDG is based upon two levels of linear problems: a global 
  system to find the solution on the boundaries of the cells, 
  followed by local systems to construct the solution inside.
  This technicality requires a careful derivation of the 
  adjoint-state method, that we address in this paper.
  We work with the acoustic wave equations in the frequency domain and 
  illustrate with a three-dimensional experiment using partial 
  reflection-data, where we further employ the features of DG-like 
  methods to efficiently handle the topography with $\mathfrak{p}$-adaptivity.
\end{abstract}

\section{Introduction}


Quantitative inverse wave problems aim to recover the
physical medium parameters that characterize the wave 
propagation from partial observations of the phenomenon.
This inverse scattering problem arises, for instance, in 
geophysics for the identification of Earth's properties, 
\cite{Lailly1983,Tarantola1984,Pratt1998,Virieux2009,Fichtner2011book},
in medical imaging or in mechanical engineering for the non-destructive 
testing, see, e.g., 
\cite{Blackledge1985,Colton1996,Scherzer2010,Cox2012,Pulkkinen2015,Ammari2015,Pham2018} 
and the references therein.

In the framework of quantitative inversion, the 
measurements of the waves (e.g., mechanical or 
electromagnetic), $\data$, are used to define a 
misfit functional $\misfit$ and the reconstruction 
of the parameters is recast as a non-linear 
minimization problem,
\begin{equation} \label{eq:min_misfit}
  \min \, \misfit(\model) \, , \qquad \text{with} \qquad
  \misfit(\model) \, = \, \dfrac{1}{2} \, \, \Big\Vert\, \forward(\model) \, - \, \data\,\Big\Vert^2 \, ,
\end{equation}
where the forward problem $\forward$ is the map from 
the model parameters $\model$ to the observable (i.e., 
the quantities measured at the position of the receivers).
That is, simulations of the wave phenomenon are compared
with the measurements to successively update the model 
parameters.
Our misfit criterion in \cref{eq:min_misfit} is the 
$L^2$ difference and least-squares problems are further 
analyzed in \cite{Chavent2010,Faucher2019IP}.
Several alternatives for the misfit have been investigated, 
we refer to, e.g., 
\cite{Luo1991,Shin2007a,Fichtner2008,Metivier2016,Yang2018,Alessandrini2019,Faucher2019FRgWI}.
Furthermore, one can incorporate a regularization term 
in \cref{eq:min_misfit} to reduce the ill-posedness, 
by the means of additional constraints, 
see, e.g., 
\cite{Rudin1992,Engl1996,Kaltenbacher2008,Kaltenbacher2018,Faucher2020EV}. 
Note that the use of a different misfit criterion 
or the incorporation of regularization terms generates
only minor modifications of the \update{methodology} we 
provide here. 

The resolution of \cref{eq:min_misfit} using deterministic 
optimization techniques makes use of algorithms in the family 
of the Newton method. 
It consists in successive updates of the model parameters 
where, at iteration $k$,
\begin{equation} \label{eq:model_update}
  \model_{k+1} \, = \, \model_{k} \, + \, \varrho_k \,\, \boldsymbol{s}_{k} \, . 
\end{equation}
In the full Newton approach, the search direction for 
the update, $\boldsymbol{s}$, depends on the gradient and 
the Hessian of the misfit functional. 
To reduce the numerical cost, alternatives that avoid the 
full Hessian computation are often employed, with Quasi 
or Truncated-Newton methods such that BFGS and 
L-BFGS \cite{Nocedal1980}, conjugate 
gradient approach \cite{Hanke1995,Metivier2013}, approximated 
pseudo-Hessian \cite{Choi2008}, or Landweber iterations. 
We refer to \cite{Nocedal2006} for a review of methods
for local optimization.
Then, the scalar step $\varrho$ is selected to obtain
an appropriate amplitude of the updates, using line-search 
algorithms, e.g., \cite{Nocedal2006,Chavent2010}.
In our implementation, we also avoid the computation of the 
Hessian such that, at each iterations of the minimization, 
one must 
\begin{enumerate}
  \item solve the forward problem using the current model parameters, 
  \item compute the gradient of the misfit functional,
\end{enumerate}
Non-linear minimization suffers from local minima, 
which cannot be avoided with the deterministic approach, 
see, e.g., \cite{Bunks1995,Sirgue2004,Faucher2019IP,Faucher2020Geo} 
in the context of seismic. 
We can mention the use of statistical-based methods but in 
the large-scale applications we have in mind, such approaches 
remain unusable at the moment.

\paragraph{Discretization of the forward problem}
The iterative minimization procedure is computationally 
intensive for large-scale applications because the 
forward problem must be solved at each iteration, 
\update{and several right-hand sides must be considered 
each time, in accordance with the number of observations.}
For instance, in seismic applications, this easily amounts 
to solving the problem for several hundreds or thousand of 
sources \update{(i.e., right-hand sides)}, at each iterations. 
For this reason, in the time-harmonic framework that we consider, 
we rely on a direct solvers (instead of iterative ones) to handle the 
linear system generated by the discretization, as it enables 
for multiple right-hand sides (rhs), that is, once the matrix 
factorization is obtained, the solution of the different rhs 
is computationally cheap. In particular, we 
use the solver \textsc{Mumps} \cite{Amestoy2001,Amestoy2006}.
On the other hand, this factorization possibly requires large 
amounts of computational memory.

Therefore, the choice of discretization methods plays a 
crucial role in the numerical efficiency.
While Finite Differences (FD) have early been employed 
for the discretization of wave problems 
(e.g., \cite{Virieux1984,Hustedt2004}), it works with a 
cartesian grid which makes it difficult to handle complex 
geometry of parameters and topography, \cite{Robertsson1996}. 
One can instead rely on methods based upon an unstructured 
mesh of the domain to have more flexibility, such that the 
Finite Element (Continuous Galerkin, CG) method 
(e.g., \cite{Chin1999,Ainsworth2006,Ern2013}), or on the 
Spectral Element Method (SEM), popularized in seismic 
applications in \cite{Komatitsch1998,Komatitsch1999}.

Then, the Discontinuous Galerkin (DG) method has been introduced,
\cite{Dumbser2006,Hesthaven2007,Brossier2010DG}.
As the name indicates, DG works with discontinuous basis 
functions, independently defined on each cell. 
It offers additional features compared to CG 
as it easily handle the $h$ and $\mathfrak{p}$-adaptivity, 
that is, the use of cells of different size and the use 
of different order of polynomial between cells, respectively, 
\cite{Hesthaven2007}. 
On the other-hand, the DG method, by working with discontinuous 
basis functions, usually results in an increase in the number of 
degrees of freedom (dof) as these are not shared between the elements.

The Hybridizable Discontinuous Galerkin method (HDG, also 
referred to as Local DG, LDG) allows to keep the features 
of the DG while avoiding the oversized linear systems. 
It shows a growing interest, e.g.,  \cite{Cockburn1998,Arnold2002,Cockburn2009,Cockburn2010,Griesmaier2011,Kirby2012,Kirby2016,Bonnasse2017,Fabien2018,Barucq2020}, 
and we refer to the introduction of \cite{Cockburn2009} 
for the background on the development of the method, and 
historical references.
Concretely, the global linear system in HDG only contains 
the dof that are on the faces of the elements, such that 
all interior ones are eliminated. 
Consequently, the method is shown to be more efficient
(less memory consumption) than CG or usual DG depending
on the order of the polynomials \cite{Kirby2012,Kirby2016,Bonnasse2017}.
Therefore, it is the perfect candidate to handle 
the large-scale inverse problem, as we can account 
for complex geometry using smaller linear systems. 
We provide in \cref{section:hdg_forward} the steps 
for the implementation: contrary to the traditional 
discretization method, it consists in two stages: 
first to compute the numerical traces (the \emph{global}
system, only composed of the dof on the faces of the elements), 
secondly, we solve the local sub-problems to build the 
volume approximation.

\paragraph{Gradient computation}
For the iterative reconstruction, the gradient of the misfit 
functional~\cref{eq:min_misfit} must be computed: 
$\nabla \misfit = D\forward^*(\forward - \data)$.
To avoid the explicit computation of the Fr\'echet derivative
$D\forward$, the \emph{adjoint-state} method has been 
designed and instead compute the action of $D\forward$, 
which is sufficient to extract the gradient, 
as we detail in \cref{section:adjoint_state_hdg}.
The method originates from the work of \cite{Lions1971} with 
early application in \cite{Chavent1974}. It is now commonly 
used and we refer to \cite{Plessix2006,Chavent2010,Pham2018,Faucher2019IP}.
With `direct' discretization method (CG, SEM, FD), the 
adjoint-state method relies on the resolution of a 
\emph{backward} problem, which reduces to the adjoint 
of the forward problem with the residuals 
(the difference $\forward - \data$) at the 
right-hand side.
Because of the two-stages in the HDG discretization
(the global and local systems), the adjoint-state
method must be appropriately derived, see 
\cref{section:adjoint_state_hdg}.

\medskip

The contribution of our work is to derive 
the methodology to solve the non-linear minimization
problem using HDG discretization, with the 
perspective of large-scale applications that are
currently intractable.
We consider the inverse acoustic wave problem 
for the identification of the medium physical 
parameters, see \cref{section:inverse_problem}.
The discretization of the wave equation using the HDG method 
is given in \cref{section:hdg_forward}, where we emphasize
the differences compared to the more traditional 
discretization methods.
In \cref{section:adjoint_state_hdg}, the adjoint-state 
method is provided in the framework of HDG, hence with 
the specificity of working with global and local sub-problems.
Eventually, we illustrate the reconstruction procedure using 
HDG with a three-dimensional experiment in \cref{section:experiments},
where we consider a medium with topography to fully use the features 
of the method, with $\mathfrak{p}$-adaptivity.
We give additional perspectives in \cref{section:perspectives}.

\section{Inverse wave problem in acoustics}
\label{section:inverse_problem}

\subsection{Time-harmonic wave propagation}

We consider the propagation of time-harmonic 
waves in acoustic medium in three 
dimensions such that $\Omega \in \R^3$, 
with boundary $\Gamma$.
We denote by $\bx$ the 3D space coordinates
$\bx=\{x,\, y,\, z\}$, and an initial (scalar)
time-harmonic source $\srcp$. 
The scalar pressure field, $\pressure: \Omega \rightarrow \Cx$, 
and the vectorial velocity, $\velocity: \Omega \rightarrow \Cx^3$, 
verify the Euler's equations,
\begin{subequations} \label{eq:euler_main}
\begin{empheq}[left={\empheqlbrace}]{align}
   -\freq \, \rho(\bx) \velocity(\bx) + \nabla \pressure(\bx)                  
               &= 0,           \quad &\text{in $\Omega$}, \label{eq:euler_main_a} \\
   -\dfrac{\freq }{\kappa(\bx)} \, \pressure(\bx) 
   +\divergence \velocity(\bx) 
               &= \srcp(\bx), \quad &\text{in $\Omega$}, \label{eq:euler_main_b}\\
  \alpha(\bx) \pressure(\bx) + \beta (\bx) \partial_{\n} \pressure(\bx) 
               &= 0,          \quad &\text{on $\Gamma$.} \label{eq:euler_main_BC}
\end{empheq} \end{subequations}
The propagation is governed by the 
physical properties of the medium: 
the density $\rho \in \R$ and the 
bulk modulus $\kappa \in \R$. 
One can also use the wave speed $c$: 
\begin{equation}
  c(\bx) = \sqrt{\kappa(\bx) \, \rho(\bx)^{-1}}.
\end{equation}
We work with a \emph{complex} frequency denoted 
$\freq$, such that 
\begin{equation} \label{eq:complex-frequency}
  \freq  \, = \, \ii \omega - \mathfrak{s} \, ,  
\end{equation}
with $\omega$ the angular frequency.
While the ``usual'' frequency-domain formulation 
would take $\mathfrak{s} = 0$, this notation is 
useful in order to work in the 
``Laplace--Fourier'' domain for the inverse problem, 
see~\cite{Shin2008,Shin2009,Faucher2017}.
It can also serve to incorporate attenuation or 
viscous behavior, as in helioseismology, e.g., \cite{Pham2019,Pham2020RRscalar}.
In the latest case, $\freq$ depends on the spatial
coordinate; this is supported in our \update{work}.
As an alternative, one can also consider \emph{complex-valued
wave speed} to account for the attenuation, see \cite{Ursin2002}.
Having complex-valued parameters does not change our \update{methodology}.

Out of generality, we have considered a Robin 
boundary condition on the boundary $\Gamma$ of 
the domain, \cref{eq:euler_main_BC} (see 
also \cref{rk:hdg-boundary-conditions}), with coefficients 
$\alpha$ and $\beta$, using $\partial_{\n}$ to 
denote the normal derivative.
It is written with respect to the pressure field 
but using~\cref{eq:euler_main}, we can equivalently 
define the Robin condition in terms of the 
velocity, or using a combination of both, such that
\begin{equation} \label{eq:euler_main_BC_velocity}
  \alpha(\bx) \pressure(\bx) 
\, + \, \freq \, \rho(\bx) \, \beta(\bx) \, \vn(\bx)
    \, = \, 0,   \qquad \text{alternative Robin condition,} 
\end{equation}
where  $\vn = \velocity \cdot \n$ indicates
the normal velocity.
In particular, in the case of \emph{absorbing 
boundary conditions} (ABC,~\cite{Engquist1977})
for acoustic media, the Robin boundary condition 
becomes
\begin{empheq}[]{align}\label{eq:abc_pressure}
 -\freq c(\bx)^{-1} \, \pressure(\bx) ~+~ \partial_{\n} \pressure(\bx)  
                   &= 0, \qquad \text{ABC in pressure,}\\
  -\big( c(\bx)\rho(\bx) \big)^{-1} \, \pressure(\bx) ~+~ 
   \vn(\bx) &= 0, \qquad \text{pressure/velocity ABC.}
   \label{eq:abc_velocity}
\end{empheq}

\begin{remark}
   One can replace the velocity field 
   in~\cref{eq:euler_main_b}
   using~\cref{eq:euler_main_a} to 
   obtain the second-order formulation 
   which only depends on
   the pressure field: 
   \begin{equation} \label{eq:euler_pressure}
       \dfrac{\freq^2}{\kappa(\bx)} \, \pressure(\bx) 
      -\divergence\Big(\dfrac{1}{\rho(\bx)}\nabla\pressure(\bx)\Big) 
     =-\freq \srcp(\bx), \quad \text{in $\Omega$.}
   \end{equation}
   In the context where the density is 
   constant, it further simplifies to the 
   Helmholtz equation. \hfill $\triangle$
\end{remark}

\subsection{Quantitative identification of the physical parameters}

In the framework of inverse problem, one wants to identify 
the model parameters that characterize the propagation, 
$\kappa$ and $\rho$, from the measurements of the waves.
In the quantitative approach, the reconstruction follows
an iterative minimization procedure, that we describe in
this section.


We first define the forward problem $\forward$ 
to give the solution of~\eqref{eq:euler_main}
at a restricted set of positions, for the model
parameters $\model = \{ \kappa, \, \rho\}$.
For a source $\srcp$ and frequency $\freq$, 
we have
\begin{equation} \label{eq:forward_problem}
\forward(\model, \,\freq, \,\srcp) = 
                            \big\{ \, \pressure (\model, \, \freq, \bx_1,              \srcp), \, \ldots , \, 
                                      \pressure (\model, \, \freq, \bx_{n_\text{rcv}}, \srcp)  \, \big\} \, ,
\end{equation}
where the $\bx_k$ are a discrete set of positions, that is, 
the forward problem gives measurements obtained from 
$n_\text{rcv}$ receivers. 
Here we have considered measurements of the pressure 
fields (commonly employed in seismic applications), 
but we can proceed similarly with the velocity, 
or with both, see, e.g., 
\cite{Alessandrini2019,Faucher2019FRgWI}.


The observed measurements are denoted by $\data$; 
these can be seen as a forward problem associated
with a target unknown model with added noise. 
The identification of parameters follows a minimization 
of a misfit functional $\misfit$ which, in the 
least-squares framework, is
\begin{equation} \label{eq:misfit}
 \misfit(\model) \, = \, \dfrac{1}{2} \, \sum_\freq \, \sum_{\update{k=1}}^{\update{N_{\mathrm{src}}}} \, 
 \big\Vert \forward(\model, \, \freq, \,\srcp_k) - \data(\freq, \, \srcp_k) \big\Vert^2_2 \,,
\end{equation}
\update{
where we consider that the data are generated from 
$N_{\mathrm{src}}$ sources, which consequently is 
the number of right-hand sides of the discretization 
problem.}
The sum over the frequencies is usually decomposed
into sub-sets, following a progression from low 
to high contents, \cite{Bunks1995,Sirgue2004,Pham2018,Faucher2020Geo}.
As mentioned in the introduction, several alternatives to the 
least-squares functional have been studied, and it can also be 
enriched with \emph{regularization} terms, see the above references.

The minimization of \cref{eq:misfit} is conducted 
following Newton-type algorithm with iterative updates
of an initial model. 
At iteration $k$, the model is updated according to 
\cref{eq:model_update}, and
we refer to \cite{Nocedal2006} for an 
extensive review of methods.
The search direction is computed from the gradient
of the misfit functional and in large-scale 
optimization, the Hessian is usually too cumbersome 
and one can use its approximation (e.g., Limited-BFGS method) 
or only the gradient (non-linear conjugate gradient methods).

For the computation of the gradient, we refer 
\cref{section:adjoint_state_hdg} where we provide
the steps for its computation with the 
adjoint-state method using HDG discretization. 
We illustrate experiments of reconstruction in
\cref{section:experiments}.

\section{Hybridizable Discontinuous Galerkin discretization}
\label{section:hdg_forward}

For numerical applications, the first step is to  
discretize the wave equation, and we follow the 
Hybridizable Discontinuous Galerkin method, that
works with the first-order problem~\cref{eq:euler_main}.
As mentioned in the introduction, this approach
has the advantage to reduce the size of the global 
linear system compared to Continuous Galerkin (depending
on the order of the polynomials), hence allowing  
to solve test-cases of larger scales.
\update{We provide here the detailed steps for the 
implementation of the method. These are necessary before we 
can develop the adjoint-state method in \cref{section:adjoint_state_hdg}.
Our methodology is established in the spirit 
of, e.g., \cite{Griesmaier2011,Bonnasse2017} (see also 
\cite{Cockburn2009,Kirby2012,Kirby2016} for elliptic problems).} 

\subsection{Notation}

\subsubsection{Domain discretization}

The domain is discretized using a non-overlapping 
partition of $\Omega$. 
The mesh of the domain is denoted $\mesh$.
It is composed of $N$ elements/cells such that
\begin{equation}
  \mesh  = \bigcup_{e=1}^\ncell \cell.
\end{equation}
The set of the $\nfaceall$ faces $\face$ is 
decomposed into the $\nfacein$ \emph{interior} 
ones (between two adjacent cells), $\facesetin$, 
and the $\nfaceout$ \emph{exterior} ones 
(between the medium and the exterior), $\facesetout$:
\begin{equation}
    \facesetall = \bigcup_{k=1}^\nfaceall \face_k \, , \qquad
    \facesetout = \facesetall \, \cap \, \Gamma \, , \qquad
    \facesetin  = \facesetall \, \setminus \, \facesetout
    \qquad \nfaceall = \nfacein + \nfaceout.
\end{equation}
In our implementation, we use simplex cells, thus triangles 
for two-dimensional domains and tetrahedra in three dimensions. 
Consequently, the face of a cell $\face$ is either a segment 
(in 2D) or a triangle (in 3D).

\subsubsection{Function and discretization spaces}

The space of discretization consists in 
piecewise polynomial of order less than 
or equal to $\mathfrak{p}$. 
In dimension three, the space of polynomial 
for simplexes is given by
\begin{equation} \label{eq:polynomial_space}
 \mathbb{P}_\mathfrak{p} = 
 \Big\{ g(\bx) = g(x,y,z) = 
 \sum_{i,j,k=0}^{\mathfrak{p}} g_{ijk} \, x^i \, y^j \, z^k \, \qquad 
                 i + j + k ~ \leq ~ \mathfrak{p} \Big\},
\end{equation}
where the $g_{ijk}$ are scalar coefficients.
Note that in DG methods, the polynomials are 
defined separately on each cell (i.e., piecewise), 
allowing discontinuities. 
We introduce the following function spaces,
associated with the polynomial on a cell $K_e$
and the cell faces $\dcell$, such that, for 
all $e$ in $1,\ldots,N$:
\begin{subequations} \label{eq:space_function_discrete}
\begin{empheq}[]{align}
 W_h \,\, &= \,\, \big\{ w_h \in L^2(\Omega), \hspace*{2.1cm}
                        w_h\mid_\cell \in \mathbb{P}_{\mathfrak{p}_e}(\cell) 
                           \, , \quad \forall \, \cell \in \mesh \big\},
                          \label{eq:space_discretization} \\
 \boldsymbol{W}_h \,\, &= \,\, \big\{ 
   \boldsymbol{w}_h \in  (L^2(\Omega))^3, \qquad \qquad 
              w_{\bullet,h}\mid_\cell \in \mathbb{P}_{\mathfrak{p}_e}(\cell)
                \, , \hspace{0.25cm} \forall \, \cell \in \mesh \big\},
                     \label{eq:space_discretization_vector} \\
  U_h \,\, &= \,\,  \big\{ u_h \in L^2(\facesetall), \hspace*{2.7cm}
                         u_h\mid_{\face_k} \in \mathbb{P}_{\mathfrak{q}_{k}}(\face_k)
                         \, , \hspace{0.6cm} \forall \, \face_k \in \facesetall \big\}. \label{eq:space_discretization_face}
\end{empheq} \end{subequations}
Therefore, we allow for a different order of polynomial
on each cell and each face. 

\subsubsection{Jump operator}

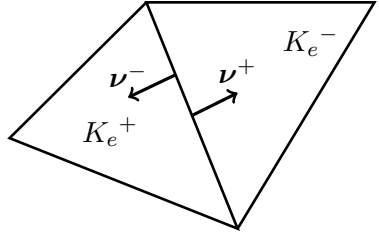
\begin{figure}[ht!] \centering
\begin{tikzpicture}[scale=1.2]
  \coordinate (t1)  at ( 0.,1.) ; \coordinate (t2) at (2.5,0.);
  \coordinate (t3)  at (1.5,2.5); \coordinate (t4) at (4. ,2.5);
  
  \coordinate (t5)  at (1.25,0.5) ; \coordinate (t6) at (0.75,1.75);
  \coordinate (t7)  at (2.0 ,1.25); \coordinate (t8) at (3.25,1.25);
  \coordinate (t9)  at (2.75,2.5);
  \draw[line width=1]  (t3) to (t1) to (t2) to (t3) to (t4) to (t2);
  \draw (1.1,0.8)  node[anchor=south] {$\cell^+$};
  \draw (3.3,1.8)  node[anchor=south] {$\cell^-$};
  \draw[line width=1.2, ->] (t7) -- (2.5,1.5);
  \draw[line width=1.2, ->] (1.82,1.70) -- (1.3,1.45);
  \draw (2.5,1.5)  node[anchor=south]  {$\n^+$};  
  \draw (1.3,1.45)  node[anchor=south] {$\n^-$};  
  \end{tikzpicture}
  \caption{Illustration of the normals 
           at the interface between two 
           triangle 
           cells.}
\label{fig:discretization:normals}
\end{figure}

We define the \emph{jump} of a quantity
between two adjacent cells, that we denote by
brackets. 
On the interface $\face$ shared by two cells (elements)
$\cell^+$ and $\cell^-$, the jump of 
$\boldsymbol{w}\cdot \n$ is
\begin{equation} \label{eq:jump}
 [\boldsymbol{w} \cdot \n]_\face 
 \, \, := \, \, \boldsymbol{w}^{\cell^+} \cdot \n_\face^+
 + \boldsymbol{w}^{\cell^-} \cdot \n_\face^-
                  = \boldsymbol{w}^{\cell^+} \cdot \n_\face^+ 
                  - \boldsymbol{w}^{\cell^-} \cdot \n_\face^+,
\end{equation}
where, by convention, $\n^{\pm}$ points outward of $\cell^{\pm}$,
as illustrated in \cref{fig:discretization:normals}.

\subsection{Local problem for the HDG discretization}

Upon assuming that the right-hand side (i.e. the source function)
$f \in L^2(\Omega)$, we can write the 
variational formulation for~\eqref{eq:euler_main_a}
and~\eqref{eq:euler_main_b}, where we 
use test functions $\phi(\bx) \in L^2(\Omega)$ and 
$\psitest(\bx) \in (L^2(\Omega))^3$.
Over each cell $\cell$ of the domain mesh, 
we have:
\begin{subequations} \begin{empheq}[left={\empheqlbrace}]{align}
  \int_\cell \Big( -\freq \, \rho \, \velocity \, \cdot \, \psitestc 
             +\nabla \pressure \,\cdot\, \psitestc \Big) ~\dd \cell &= 0,  \\
  \int_\cell \Big( -\freq \, \kappa^{-1} \, \pressure \, \phic
             +(\divergence \velocity) \, \phic \, \Big) ~\dd \cell
                         &= \int_\cell \srcp \, \phic ~\dd \cell,
\end{empheq} \end{subequations}
where $\overline{\phantom{\psi}}$ denotes 
the conjugation and we omit the space dependency for clarity.
Integration by part gives 
\begin{subequations} \label{eq:main_variational_ipp}
\begin{empheq}[left={\empheqlbrace}]{align}
\int_\cell \Big( -\freq \, \rho \, \velocity \,\cdot\, \psitestc 
           - \pressure  \, (\divergence \psitestc) \Big)~\dd \cell \, 
 + \, \int_\dcell \pressure \,\, \psitestc \cdot \n ~\dd \dcell &= 0 \, ,  \\
  \int_\cell \Big( -\freq \, \kappa^{-1} \, \pressure \, \phic
             -\velocity \cdot \nabla \phic 
             -\srcp \, \phic \Big) ~\dd \cell
 +\, \int_\dcell \phic \,\, \velocity \cdot \n ~ \dd\dcell  &= 0 \, .
\end{empheq} \end{subequations}
Note that the integral over the cell boundary
can be substituted by a sum over the faces of
the cell such that, e.g., 
\begin{equation} \label{eq:surface-int_to_sum-face}
  \int_\dcell \pressure \,\, \psitestc \cdot \n ~\dd \dcell \, = \,
  \sum_{\face \in \dcell} \int_\face \pressure \,\, \psitestc \cdot \n_\face ~\dd \face \, .
\end{equation}

\medskip

The solutions $\pressure$ and $\velocity$ are 
approximated by polynomials of order $\mathfrak{p}_e$ 
on the cell, respectively with the discretization 
variables 
\begin{equation}
  \pressure_h \in W_h 
  \qquad \text{~and~} \qquad 
  \velocity_h \in \boldsymbol{W}_h \, ,
\end{equation}
and we refer to $\widehat{\pressure}_h \in U_h$ 
for the numerical trace of $\pressure$.

By introducing the discretization variables 
in \cref{eq:main_variational_ipp}, we obtain,
\begin{subequations}  \label{eq:variational_approx-var}
\begin{empheq}[left={\empheqlbrace}]{align} 
& \int_\cell\Big( -\freq  \rho  \velocity_h^{(e)} \cdot \psitestc 
           - \pressure_h^{(e)}   (\divergence \psitestc) \Big) \, \dd \cell  
 + \int_\dcell \widehat{\pressure}_h^{(e)}  \psitestc 
              \cdot \n \, \dd \dcell = 0, \label{eq:variational_approx-var_1} \\
& \int_\cell \Big( -\freq  \kappa^{-1}  \pressure_h^{(e)}  \phic
             -\velocity_h^{(e)} \cdot \nabla \phic 
             -\srcp  \phic \Big) \, \dd \cell
 +  \int_\dcell \hspace*{-4mm} \phic  \widehat{\velocity}_h^{(e)} \cdot \n \, \dd \dcell = 0,
 \label{eq:variational_approx-var_2}
\end{empheq} \end{subequations}
where $\pressure_h^{(e)} \in \mathbb{P}_{\mathfrak{p}_e} = \pressure_h \mid_{\cell}$,
similarly for $\velocity_h^{(e)}$, according to~\cref{eq:space_function_discrete}.
The essence of the HDG method is to formulate the 
\emph{numerical flux} $\widehat{\velocity}_h^{(e)}$ 
such that (\cite{Cockburn2008,Cockburn2009,Cockburn2010})
\begin{equation} \label{eq:approx_trace}
  \widehat{\velocity}_h^{(e)} \, = \, 
          \velocity_h^{(e)} 
  \, + \, \tau \Big(\pressure_h^{(e)}  - \lambda_h^{(e)}\Big) \, \n^{(e)} \, ,
  \qquad \text{with } \quad \lambda_h^{(e)} \, = \, \widehat{\pressure}_h^{(e)} \, ,  
\end{equation}
where $\n^{(e)}$ denotes the normal on
the boundary of $\cell$.
Here, $\tau$ is a penalization/stabilization parameter defined 
on each face and allowed to have distinct values for the two 
sides of the faces, that is, as illustrated 
in \cref{fig:discretization:normals}, on a face $\face$ we 
distinguish between $\tau_\face^{+}$ and $\tau_\face^{-}$.
It is required that $\tau > 0$ for uniqueness, see \cite{Arnold2002,Cockburn2009,Nguyen2009}.
Practically, following the work of \cite{Nguyen2009}
for diffusion problems, we use
\begin{equation}
 \tau_\face^{\pm} = \rho^{-1}(\cell^{\pm}) \, , \quad
 \text{for a piecewise constant (per cell) density model.}
\end{equation}
On the one hand, $\pressure_h$ and $\velocity_h$ 
are piecewise polynomial on the cells, thus allowing 
discontinuities between two adjacent cells 
(\cref{fig:discretization:normals}); on the other 
hand, $\lambda_h$ is defined on the `skeleton' of the 
mesh (i.e., on the faces of the cells) and allows 
discontinuity only at the nodes in two dimensions, 
and at the edges in three dimensions, see \cref{fig:dof}.

\medskip

In the following, we omit the cell exponent $^{(e)}$ for the 
sake of clarity. 
\update{We replace the numerical flux in \cref{eq:variational_approx-var_2}
by \cref{eq:approx_trace}, and using the revert integration by part, 
we obtain,
\begin{equation} \label{eq:variational_discretization_B}
  \int_\cell \Big( -\freq \, \kappa^{-1} \, \pressure_h \, \phic
           +(\divergence \velocity_h)\, \phic - \srcp \, \phic \Big)~\dd \cell 
+ \int_\dcell \tau (\pressure_h  - \lambda_h) \, \, \phic ~\dd\dcell = 0.
\end{equation}}
\subsection{Continuity condition for the HDG discretization}

To complete the problem, the continuity of the discretized 
velocity field is enforced and we impose that, on all of the 
cell faces, the jump \cref{eq:jump} equates zero. 
For any test function $\xi \in U_h$, we write
\begin{equation} \label{eq:variational_discretization_additional} 
  \forall {\face \in \facesetall}, \quad \int_{\face}  
  \big[ \widehat{\velocity}_h \cdot \n_\face \big] \, \xic ~\dd \face = 0 \, ,
\end{equation}
with $\n_\face$ the normal on the face.
For an inner boundary, we insert~\cref{eq:approx_trace}
and use the jump definition~\cref{eq:jump}: 
\begin{empheq}[]{align} \label{eq:lambda_boundary-element}
\begin{split}
   \forall {\face \in \facesetin}, \qquad 
&  \int_{\face} \big[ \widehat{\velocity}_h \cdot \n_\face \big] \, \xic ~\dd \face 
  =\int_{\face} \, \Big( \, \widehat{\velocity}_h^{(\face^+)} \cdot \n_\face^+
                       \, \, + \, \, \widehat{\velocity}_h^{(\face^-)} \cdot \n_\face^- \, \Big) \xic \, \dd \face \\
& \hspace*{-1.5cm} = \int_{\face} \, \Big( \, 
   \velocity_h^{+} \cdot \n_\face^+ \, + \, \tau^{+} (\pressure_h^{+} - \lambda_h) 
  \, + \,
   \velocity_h^{-} \cdot \n_\face^- \, + \, \tau^{-} (\pressure_h^{-} - \lambda_h) 
  \, \Big) \, \xic ~\dd \face \, ,
\end{split} \end{empheq}
where the exponent $^{\pm}$ indicates
the adjacent cells, cf. \cref{fig:discretization:normals}.
On a face that belongs to the outer boundary, only one 
side remains, and we have 
\begin{equation} \label{eq:boundary-domain}
\forall {\face \in \facesetout}, \qquad 
\widehat{\velocity}_h \mid_{\face \in \facesetout}
 \, = \, \velocity_h 
        \, + \, \tau \big(\pressure_h  - \lambda_h \big) \, \n_{\face} \, .
\end{equation}
Then, using the specified Robin boundary condition, 
$\widehat{\velocity}_h$ must also 
verifies~\cref{eq:euler_main_BC_velocity},
such that, 
\begin{equation} \label{eq:lambda_boundary-domain}
\forall {\face \in \facesetout}, \qquad 
  \widehat{\velocity}_h \, \cdot \, \n_\face
  \overset{\cref{eq:euler_main_BC_velocity}}{=}
 -\dfrac{\alpha}{\freq \, \rho \, \beta}  \, \lambda_h
  \overset{\cref{eq:boundary-domain}}{=}
   \velocity_h \, \cdot \, \n_{\face} \, + \, \tau \big(\pressure_h  - \lambda_h \big) \, ,
\end{equation}
where we also replace the trace of $\pressure$
by $\lambda_h$ in the Robin condition.
\medskip

Eventually, we summarize the complete HDG problem.
Find $(\lambda_h, \, \pressure_h, \, \velocity_h)$ that 
solve 
\cref{eq:variational_approx-var_1} (inserting $\lambda_h$), 
\cref{eq:variational_discretization_B,eq:lambda_boundary-element,eq:lambda_boundary-domain}, 
that is, using~\cref{eq:surface-int_to_sum-face}:
\begin{equation}\label{eq:variational_discretization_full}
\begin{aligned}
& \int_\cell \Big( -\freq  \rho  \velocity_h \cdot \psitestc 
           - \pressure_h   (\divergence \psitestc) \Big)~\dd \cell  
 + \sum_{\face \in \dcell} \int_\face  \lambda_h  \psitestc 
              \cdot \n_\face ~\dd \face = 0  , \hspace*{0.20cm} \forall \cell \in \mesh  , \\[0.5em]
& \int_\cell \Big(
           -\freq \kappa^{-1} \pressure_h \phic
           +(\divergence \velocity_h)\phic - \srcp \phic \Big) \dd \cell  
+ \sum_{\face \in \dcell} \int_\face  \tau_\face (\pressure_h  - \lambda_h ) 
              \phic  \dd\face =  0 ,  
              \forall \cell, \\[0.5em]
& \int_{\face}  \Big(  
   \velocity_h^{+} \cdot \n_\face^+ + \tau_\face^{+} (\pressure_h^{+} - \lambda_h) 
    + 
   \velocity_h^{-} \cdot \n_\face^- + \tau_\face^{-} (\pressure_h^{-} - \lambda_h) 
   \Big)  \xic ~\dd \face  = 0  , \hspace*{2mm} 
          \forall {\face \in \facesetin}, \\[0.5em]
& \int_{\face} 
 \bigg( \velocity_h  \cdot  \n_{\face}  +  
 \tau_\face  \pressure_h      
         +  
        \Big( \dfrac{\alpha}{\freq  \rho  \beta}  -  \tau_\face \Big) \lambda_h
  \bigg)  \xic ~\dd \face ~ = 0 ,
             \hspace*{3.20cm}  \forall {\face \in \facesetout}.
\end{aligned}
\end{equation}
The two first equations represent \emph{local} problems, 
on each cell, while the last two equations give the 
conditions on the boundaries of the cells.

\subsection{Linear systems}

\begin{figure}[ht!] \centering
  \subcaptionbox{Degrees of freedom for Continuous 
                 Galerkin discretization \update{(16 total)}.
                 \label{fig:dof:fe}}
                {\begin{tikzpicture}[scale=1.05]
  \pgfmathsetmacro{\xa}{0.0}   \pgfmathsetmacro{\ya}{1.0}
  \pgfmathsetmacro{\xb}{2.5}   \pgfmathsetmacro{\yb}{0.0}
  \pgfmathsetmacro{\xc}{1.5}   \pgfmathsetmacro{\yc}{2.5}
  \pgfmathsetmacro{\xd}{4.0}   \pgfmathsetmacro{\yd}{2.5}
  \coordinate (t1) at (\xa,\ya) ;
  \coordinate (t2) at (\xb,\yb) ;
  \coordinate (t3) at (\xc,\yc) ; 
  \coordinate (t4) at (\xd,\yd) ; 
  \draw[line width=1, mark size=1.25pt]  (t3) to (t1) to (t2) to (t3) to (t4) to (t2);

  \node[line width=1, mark size=1.25pt] at (t1) {\pgfuseplotmark{*}};
  \node[line width=1, mark size=1.25pt] at (t2) {\pgfuseplotmark{*}};
  \node[line width=1, mark size=1.25pt] at (t3) {\pgfuseplotmark{*}};
  \node[line width=1, mark size=1.25pt] at (t4) {\pgfuseplotmark{*}};

  \pgfmathsetmacro{\ox}{1./3.*(\xa+\xb+\xc)}
  \pgfmathsetmacro{\oy}{1./3.*(\ya+\yb+\yc)}
  \node[line width=1, mark size=1.25pt] at (\ox,\oy) {\pgfuseplotmark{*}};  
  \pgfmathsetmacro{\ox}{1./3.*(\xd+\xb+\xc)}
  \pgfmathsetmacro{\oy}{1./3.*(\yd+\yb+\yc)}
  \node[line width=1, mark size=1.25pt] at (\ox,\oy) {\pgfuseplotmark{*}};

  \pgfmathsetmacro{\ox}{\xa}   \pgfmathsetmacro{\oy}{\ya}
  \pgfmathsetmacro{\nx}{\xb}   \pgfmathsetmacro{\ny}{\yb}
  \pgfmathsetmacro{\sizeedgeX}{\nx-\ox} \pgfmathsetmacro{\sizeedgeY}{\ny-\oy}
  \pgfmathsetmacro{\myx} {\ox + \sizeedgeX/3.} 
  \pgfmathsetmacro{\myy} {\oy + \sizeedgeY/3.}
  \node[line width=1, mark size=1.25pt] at (\myx,\myy) {\pgfuseplotmark{*}};  
  \pgfmathsetmacro{\myx} {\ox + + 2*\sizeedgeX/3.}
  \pgfmathsetmacro{\myy} {\oy + + 2*\sizeedgeY/3.}
  \node[line width=1, mark size=1.25pt] at (\myx,\myy) {\pgfuseplotmark{*}};  
  
  \pgfmathsetmacro{\ox}{\xa}   \pgfmathsetmacro{\oy}{\ya}
  \pgfmathsetmacro{\nx}{\xc}   \pgfmathsetmacro{\ny}{\yc}
  \pgfmathsetmacro{\sizeedgeX}{\nx-\ox} \pgfmathsetmacro{\sizeedgeY}{\ny-\oy}
  \pgfmathsetmacro{\myx} {\ox + \sizeedgeX/3.} 
  \pgfmathsetmacro{\myy} {\oy + \sizeedgeY/3.}
  \node[line width=1, mark size=1.25pt] at (\myx,\myy) {\pgfuseplotmark{*}};  
  \pgfmathsetmacro{\myx} {\ox + + 2*\sizeedgeX/3.}
  \pgfmathsetmacro{\myy} {\oy + + 2*\sizeedgeY/3.}
  \node[line width=1, mark size=1.25pt] at (\myx,\myy) {\pgfuseplotmark{*}};  
  
  \pgfmathsetmacro{\ox}{\xb}   \pgfmathsetmacro{\oy}{\yb}
  \pgfmathsetmacro{\nx}{\xc}   \pgfmathsetmacro{\ny}{\yc}
  \pgfmathsetmacro{\sizeedgeX}{\nx-\ox} \pgfmathsetmacro{\sizeedgeY}{\ny-\oy}
  \pgfmathsetmacro{\myx} {\ox + \sizeedgeX/3.} 
  \pgfmathsetmacro{\myy} {\oy + \sizeedgeY/3.}
  \node[line width=1, mark size=1.25pt] at (\myx,\myy) {\pgfuseplotmark{*}};  
  \pgfmathsetmacro{\myx} {\ox + + 2*\sizeedgeX/3.}
  \pgfmathsetmacro{\myy} {\oy + + 2*\sizeedgeY/3.}
  \node[line width=1, mark size=1.25pt] at (\myx,\myy) {\pgfuseplotmark{*}};  
  
  \pgfmathsetmacro{\ox}{\xb}   \pgfmathsetmacro{\oy}{\yb}
  \pgfmathsetmacro{\nx}{\xd}   \pgfmathsetmacro{\ny}{\yd}
  \pgfmathsetmacro{\sizeedgeX}{\nx-\ox} \pgfmathsetmacro{\sizeedgeY}{\ny-\oy}
  \pgfmathsetmacro{\myx} {\ox + \sizeedgeX/3.} 
  \pgfmathsetmacro{\myy} {\oy + \sizeedgeY/3.}
  \node[line width=1, mark size=1.25pt] at (\myx,\myy) {\pgfuseplotmark{*}};  
  \pgfmathsetmacro{\myx} {\ox + + 2*\sizeedgeX/3.}
  \pgfmathsetmacro{\myy} {\oy + + 2*\sizeedgeY/3.}
  \node[line width=1, mark size=1.25pt] at (\myx,\myy) {\pgfuseplotmark{*}};  
   
  \pgfmathsetmacro{\ox}{\xc}   \pgfmathsetmacro{\oy}{\yc}
  \pgfmathsetmacro{\nx}{\xd}   \pgfmathsetmacro{\ny}{\yd}
  \pgfmathsetmacro{\sizeedgeX}{\nx-\ox} \pgfmathsetmacro{\sizeedgeY}{\ny-\oy}
  \pgfmathsetmacro{\myx} {\ox + \sizeedgeX/3.} 
  \pgfmathsetmacro{\myy} {\oy + \sizeedgeY/3.}
  \node[line width=1, mark size=1.25pt] at (\myx,\myy) {\pgfuseplotmark{*}};  
  \pgfmathsetmacro{\myx} {\ox + + 2*\sizeedgeX/3.}
  \pgfmathsetmacro{\myy} {\oy + + 2*\sizeedgeY/3.}
  \node[line width=1, mark size=1.25pt] at (\myx,\myy) {\pgfuseplotmark{*}};  

\end{tikzpicture}
}  \hfill 
  \subcaptionbox{Degrees of freedom for $\pressure_h$ and $\velocity_h$ 
                 (i.e., Internal Penalty Discontinuous Galerkin discretization)
                 \update{(20 total)}.
                 \label{fig:dof:dg}}
                {\begin{tikzpicture}[scale=1.05]
  \pgfmathsetmacro{\xa}{0.0}   \pgfmathsetmacro{\ya}{1.0}
  \pgfmathsetmacro{\xb}{2.5}   \pgfmathsetmacro{\yb}{0.0}
  \pgfmathsetmacro{\xc}{1.5}   \pgfmathsetmacro{\yc}{2.5}
  \pgfmathsetmacro{\xd}{4.0}   \pgfmathsetmacro{\yd}{2.5}
  \coordinate (t1) at (\xa,\ya) ;
  \coordinate (t2) at (\xb,\yb) ;
  \coordinate (t3) at (\xc,\yc) ; 
  \coordinate (t4) at (\xd,\yd) ; 
  \draw[line width=1, mark size=1.25pt]  (t3) to (t1) to (t2) to (t3) to (t4) to (t2);
  \pgfmathsetmacro{\ox}{1./3.*(\xa+\xb+\xc)}
  \pgfmathsetmacro{\oy}{1./3.*(\ya+\yb+\yc)}
  \node[line width=1, mark size=1.25pt] at (\ox,\oy) {\pgfuseplotmark{*}};  
  \pgfmathsetmacro{\ox}{1./3.*(\xd+\xb+\xc)}
  \pgfmathsetmacro{\oy}{1./3.*(\yd+\yb+\yc)}
  \node[line width=1, mark size=1.25pt] at (\ox,\oy) {\pgfuseplotmark{*}};

  \pgfmathsetmacro{\xa}{0.25}  \pgfmathsetmacro{\ya}{1.05}
  \pgfmathsetmacro{\xb}{2.28}  \pgfmathsetmacro{\yb}{0.22}
  \pgfmathsetmacro{\xc}{1.45}  \pgfmathsetmacro{\yc}{2.25}
  \coordinate (t1) at (\xa,\ya) ;
  \coordinate (t2) at (\xb,\yb) ;
  \coordinate (t3) at (\xc,\yc) ; 
  \node[line width=1, mark size=1.25pt] at (t1) {\pgfuseplotmark{*}};
  \node[line width=1, mark size=1.25pt] at (t2) {\pgfuseplotmark{*}};
  \node[line width=1, mark size=1.25pt] at (t3) {\pgfuseplotmark{*}};  
  \pgfmathsetmacro{\ox}{\xa}   \pgfmathsetmacro{\oy}{\ya}
  \pgfmathsetmacro{\nx}{\xb}   \pgfmathsetmacro{\ny}{\yb}
  \pgfmathsetmacro{\sizeedgeX}{\nx-\ox} \pgfmathsetmacro{\sizeedgeY}{\ny-\oy}
  \pgfmathsetmacro{\myx} {\ox + \sizeedgeX/3.} 
  \pgfmathsetmacro{\myy} {\oy + \sizeedgeY/3.}
  \node[line width=1, mark size=1.25pt] at (\myx,\myy) {\pgfuseplotmark{*}};  
  \pgfmathsetmacro{\myx} {\ox + + 2*\sizeedgeX/3.}
  \pgfmathsetmacro{\myy} {\oy + + 2*\sizeedgeY/3.}
  \node[line width=1, mark size=1.25pt] at (\myx,\myy) {\pgfuseplotmark{*}};  
  \pgfmathsetmacro{\ox}{\xa}   \pgfmathsetmacro{\oy}{\ya}
  \pgfmathsetmacro{\nx}{\xc}   \pgfmathsetmacro{\ny}{\yc}
  \pgfmathsetmacro{\sizeedgeX}{\nx-\ox} \pgfmathsetmacro{\sizeedgeY}{\ny-\oy}
  \pgfmathsetmacro{\myx} {\ox + \sizeedgeX/3.} 
  \pgfmathsetmacro{\myy} {\oy + \sizeedgeY/3.}
  \node[line width=1, mark size=1.25pt] at (\myx,\myy) {\pgfuseplotmark{*}};  
  \pgfmathsetmacro{\myx} {\ox + + 2*\sizeedgeX/3.}
  \pgfmathsetmacro{\myy} {\oy + + 2*\sizeedgeY/3.}
  \node[line width=1, mark size=1.25pt] at (\myx,\myy) {\pgfuseplotmark{*}};  
  \pgfmathsetmacro{\ox}{\xb}   \pgfmathsetmacro{\oy}{\yb}
  \pgfmathsetmacro{\nx}{\xc}   \pgfmathsetmacro{\ny}{\yc}
  \pgfmathsetmacro{\sizeedgeX}{\nx-\ox} \pgfmathsetmacro{\sizeedgeY}{\ny-\oy}
  \pgfmathsetmacro{\myx} {\ox + \sizeedgeX/3.} 
  \pgfmathsetmacro{\myy} {\oy + \sizeedgeY/3.}
  \node[line width=1, mark size=1.25pt] at (\myx,\myy) {\pgfuseplotmark{*}};  
  \pgfmathsetmacro{\myx} {\ox + + 2*\sizeedgeX/3.}
  \pgfmathsetmacro{\myy} {\oy + + 2*\sizeedgeY/3.}
  \node[line width=1, mark size=1.25pt] at (\myx,\myy) {\pgfuseplotmark{*}};

  \pgfmathsetmacro{\xa}{3.78}  \pgfmathsetmacro{\ya}{2.36}
  \pgfmathsetmacro{\xb}{2.52}  \pgfmathsetmacro{\yb}{0.28}
  \pgfmathsetmacro{\xc}{1.75}  \pgfmathsetmacro{\yc}{2.35}
  \coordinate (t1) at (\xa,\ya) ;
  \coordinate (t2) at (\xb,\yb) ;
  \coordinate (t3) at (\xc,\yc) ; 
  \node[line width=1, mark size=1.25pt] at (t1) {\pgfuseplotmark{*}};
  \node[line width=1, mark size=1.25pt] at (t2) {\pgfuseplotmark{*}};
  \node[line width=1, mark size=1.25pt] at (t3) {\pgfuseplotmark{*}};  
  \pgfmathsetmacro{\ox}{\xa}   \pgfmathsetmacro{\oy}{\ya}
  \pgfmathsetmacro{\nx}{\xb}   \pgfmathsetmacro{\ny}{\yb}
  \pgfmathsetmacro{\sizeedgeX}{\nx-\ox} \pgfmathsetmacro{\sizeedgeY}{\ny-\oy}
  \pgfmathsetmacro{\myx} {\ox + \sizeedgeX/3.} 
  \pgfmathsetmacro{\myy} {\oy + \sizeedgeY/3.}
  \node[line width=1, mark size=1.25pt] at (\myx,\myy) {\pgfuseplotmark{*}};  
  \pgfmathsetmacro{\myx} {\ox + + 2*\sizeedgeX/3.}
  \pgfmathsetmacro{\myy} {\oy + + 2*\sizeedgeY/3.}
  \node[line width=1, mark size=1.25pt] at (\myx,\myy) {\pgfuseplotmark{*}};  
  \pgfmathsetmacro{\ox}{\xa}   \pgfmathsetmacro{\oy}{\ya}
  \pgfmathsetmacro{\nx}{\xc}   \pgfmathsetmacro{\ny}{\yc}
  \pgfmathsetmacro{\sizeedgeX}{\nx-\ox} \pgfmathsetmacro{\sizeedgeY}{\ny-\oy}
  \pgfmathsetmacro{\myx} {\ox + \sizeedgeX/3.} 
  \pgfmathsetmacro{\myy} {\oy + \sizeedgeY/3.}
  \node[line width=1, mark size=1.25pt] at (\myx,\myy) {\pgfuseplotmark{*}};  
  \pgfmathsetmacro{\myx} {\ox + + 2*\sizeedgeX/3.}
  \pgfmathsetmacro{\myy} {\oy + + 2*\sizeedgeY/3.}
  \node[line width=1, mark size=1.25pt] at (\myx,\myy) {\pgfuseplotmark{*}};  
  \pgfmathsetmacro{\ox}{\xb}   \pgfmathsetmacro{\oy}{\yb}
  \pgfmathsetmacro{\nx}{\xc}   \pgfmathsetmacro{\ny}{\yc}
  \pgfmathsetmacro{\sizeedgeX}{\nx-\ox} \pgfmathsetmacro{\sizeedgeY}{\ny-\oy}
  \pgfmathsetmacro{\myx} {\ox + \sizeedgeX/3.} 
  \pgfmathsetmacro{\myy} {\oy + \sizeedgeY/3.}
  \node[line width=1, mark size=1.25pt] at (\myx,\myy) {\pgfuseplotmark{*}};  
  \pgfmathsetmacro{\myx} {\ox + + 2*\sizeedgeX/3.}
  \pgfmathsetmacro{\myy} {\oy + + 2*\sizeedgeY/3.}
  \node[line width=1, mark size=1.25pt] at (\myx,\myy) {\pgfuseplotmark{*}}; 
\end{tikzpicture}
}\hfill 
  \subcaptionbox{Degrees of freedom for $\lambda_h$ which give the global 
                 matrix size in the HDG discretization
                 \update{(20 total)}. \label{fig:dof:hdg}}
                 {\begin{tikzpicture}[scale=1.05]
  \pgfmathsetmacro{\step}{0.1}
  \pgfmathsetmacro{\xa}{0.0}   \pgfmathsetmacro{\ya}{1.0}
  \pgfmathsetmacro{\xb}{2.5}   \pgfmathsetmacro{\yb}{0.0}
  \pgfmathsetmacro{\xc}{1.5}   \pgfmathsetmacro{\yc}{2.5}
  \pgfmathsetmacro{\xd}{4.0}   \pgfmathsetmacro{\yd}{2.5}
  \coordinate (t1) at (\xa,\ya) ;
  \coordinate (t2) at (\xb,\yb) ;
  \coordinate (t3) at (\xc,\yc) ; 
  \coordinate (t4) at (\xd,\yd) ; 
  \coordinate (t1a) at (\xa,\ya+\step)  ;
  \coordinate (t1b) at (\xa,\ya-\step)  ;
  \coordinate (t2a) at (\xb-\step,\yb)  ;
  \coordinate (t2b) at (\xb+\step,\yb)  ;
  \coordinate (t2c) at (\xb+3*\step,\yb);
  \coordinate (t3a) at (\xc-\step,\yc)  ; 
  \coordinate (t3b) at (\xc+\step,\yc)  ; 
  \coordinate (t3c) at (\xc+3*\step,\yc); 
  \coordinate (t4a) at (\xd-\step,\yd)  ; 
  \coordinate (t4b) at (\xd,\yd-2*\step)  ; 
  
  \draw[line width=1, mark size=1.25pt]  (t3a) to (t1a) ;
  \draw[line width=1, mark size=1.25pt]  (t1b) to (t2a) ;
  \draw[line width=1, mark size=1.25pt]  (t2b) to (t3b) ;
  \draw[line width=1, mark size=1.25pt]  (t3c) to (t4a) ;
  \draw[line width=1, mark size=1.25pt]  (t4b) to (t2c) ;

  \pgfmathsetmacro{\ox}{\xa}       \pgfmathsetmacro{\oy}{\ya-\step}
  \pgfmathsetmacro{\nx}{\xb-\step} \pgfmathsetmacro{\ny}{\yb}
  \pgfmathsetmacro{\sizeedgeX}{\nx-\ox} \pgfmathsetmacro{\sizeedgeY}{\ny-\oy}
  \pgfmathsetmacro{\myx} {\ox + \sizeedgeX/3.} 
  \pgfmathsetmacro{\myy} {\oy + \sizeedgeY/3.}
  \node[line width=1, mark size=1.25pt] at (\myx,\myy) {\pgfuseplotmark{*}};  
  \pgfmathsetmacro{\myx} {\ox + + 2*\sizeedgeX/3.}
  \pgfmathsetmacro{\myy} {\oy + + 2*\sizeedgeY/3.}
  \node[line width=1, mark size=1.25pt] at (\myx,\myy) {\pgfuseplotmark{*}};  
  \pgfmathsetmacro{\myx} {\ox}
  \pgfmathsetmacro{\myy} {\oy}
  \node[line width=1, mark size=1.25pt] at (\myx,\myy) {\pgfuseplotmark{*}};  
  \pgfmathsetmacro{\myx} {\nx}
  \pgfmathsetmacro{\myy} {\ny}
  \node[line width=1, mark size=1.25pt] at (\myx,\myy) {\pgfuseplotmark{*}};  

  \pgfmathsetmacro{\ox}{\xa}       \pgfmathsetmacro{\oy}{\ya+\step}
  \pgfmathsetmacro{\nx}{\xc-\step} \pgfmathsetmacro{\ny}{\yc}
  \pgfmathsetmacro{\sizeedgeX}{\nx-\ox} \pgfmathsetmacro{\sizeedgeY}{\ny-\oy}
  \pgfmathsetmacro{\myx} {\ox + \sizeedgeX/3.} 
  \pgfmathsetmacro{\myy} {\oy + \sizeedgeY/3.}
  \node[line width=1, mark size=1.25pt] at (\myx,\myy) {\pgfuseplotmark{*}};  
  \pgfmathsetmacro{\myx} {\ox + + 2*\sizeedgeX/3.}
  \pgfmathsetmacro{\myy} {\oy + + 2*\sizeedgeY/3.}
  \node[line width=1, mark size=1.25pt] at (\myx,\myy) {\pgfuseplotmark{*}};  
  \pgfmathsetmacro{\myx} {\ox}
  \pgfmathsetmacro{\myy} {\oy}
  \node[line width=1, mark size=1.25pt] at (\myx,\myy) {\pgfuseplotmark{*}};  
  \pgfmathsetmacro{\myx} {\nx}
  \pgfmathsetmacro{\myy} {\ny}
  \node[line width=1, mark size=1.25pt] at (\myx,\myy) {\pgfuseplotmark{*}};  

  \pgfmathsetmacro{\ox}{\xb+\step} \pgfmathsetmacro{\oy}{\yb}
  \pgfmathsetmacro{\nx}{\xc+\step} \pgfmathsetmacro{\ny}{\yc}
  \pgfmathsetmacro{\sizeedgeX}{\nx-\ox} \pgfmathsetmacro{\sizeedgeY}{\ny-\oy}
  \pgfmathsetmacro{\myx} {\ox + \sizeedgeX/3.} 
  \pgfmathsetmacro{\myy} {\oy + \sizeedgeY/3.}
  \node[line width=1, mark size=1.25pt] at (\myx,\myy) {\pgfuseplotmark{*}};  
  \pgfmathsetmacro{\myx} {\ox + + 2*\sizeedgeX/3.}
  \pgfmathsetmacro{\myy} {\oy + + 2*\sizeedgeY/3.}
  \node[line width=1, mark size=1.25pt] at (\myx,\myy) {\pgfuseplotmark{*}};  
  \pgfmathsetmacro{\myx} {\ox}
  \pgfmathsetmacro{\myy} {\oy}
  \node[line width=1, mark size=1.25pt] at (\myx,\myy) {\pgfuseplotmark{*}};  
  \pgfmathsetmacro{\myx} {\nx}
  \pgfmathsetmacro{\myy} {\ny}
  \node[line width=1, mark size=1.25pt] at (\myx,\myy) {\pgfuseplotmark{*}};  

  \pgfmathsetmacro{\ox}{\xc+3*\step} \pgfmathsetmacro{\oy}{\yc}
  \pgfmathsetmacro{\nx}{\xd-\step}   \pgfmathsetmacro{\ny}{\yd}
  \pgfmathsetmacro{\sizeedgeX}{\nx-\ox} \pgfmathsetmacro{\sizeedgeY}{\ny-\oy}
  \pgfmathsetmacro{\myx} {\ox + \sizeedgeX/3.} 
  \pgfmathsetmacro{\myy} {\oy + \sizeedgeY/3.}
  \node[line width=1, mark size=1.25pt] at (\myx,\myy) {\pgfuseplotmark{*}};  
  \pgfmathsetmacro{\myx} {\ox + + 2*\sizeedgeX/3.}
  \pgfmathsetmacro{\myy} {\oy + + 2*\sizeedgeY/3.}
  \node[line width=1, mark size=1.25pt] at (\myx,\myy) {\pgfuseplotmark{*}};  
  \pgfmathsetmacro{\myx} {\ox}
  \pgfmathsetmacro{\myy} {\oy}
  \node[line width=1, mark size=1.25pt] at (\myx,\myy) {\pgfuseplotmark{*}};  
  \pgfmathsetmacro{\myx} {\nx}
  \pgfmathsetmacro{\myy} {\ny}
  \node[line width=1, mark size=1.25pt] at (\myx,\myy) {\pgfuseplotmark{*}};  

  \pgfmathsetmacro{\ox}{\xd}         \pgfmathsetmacro{\oy}{\yd-2*\step}
  \pgfmathsetmacro{\nx}{\xb+3*\step} \pgfmathsetmacro{\ny}{\yb}
  \pgfmathsetmacro{\sizeedgeX}{\nx-\ox} \pgfmathsetmacro{\sizeedgeY}{\ny-\oy}
  \pgfmathsetmacro{\myx} {\ox + \sizeedgeX/3.} 
  \pgfmathsetmacro{\myy} {\oy + \sizeedgeY/3.}
  \node[line width=1, mark size=1.25pt] at (\myx,\myy) {\pgfuseplotmark{*}};  
  \pgfmathsetmacro{\myx} {\ox + + 2*\sizeedgeX/3.}
  \pgfmathsetmacro{\myy} {\oy + + 2*\sizeedgeY/3.}
  \node[line width=1, mark size=1.25pt] at (\myx,\myy) {\pgfuseplotmark{*}};  
  \pgfmathsetmacro{\myx} {\ox}
  \pgfmathsetmacro{\myy} {\oy}
  \node[line width=1, mark size=1.25pt] at (\myx,\myy) {\pgfuseplotmark{*}};  
  \pgfmathsetmacro{\myx} {\nx}
  \pgfmathsetmacro{\myy} {\ny}
  \node[line width=1, mark size=1.25pt] at (\myx,\myy) {\pgfuseplotmark{*}};  
\end{tikzpicture}
}
  \caption{Illustration of the degrees of freedom in two dimensions
           \update{for polynomial order 3. By avoiding the inner degrees
           of freedom, HDG is more efficient (i.e., reduces their number)
           at high order.}}
  \label{fig:dof}
\end{figure}
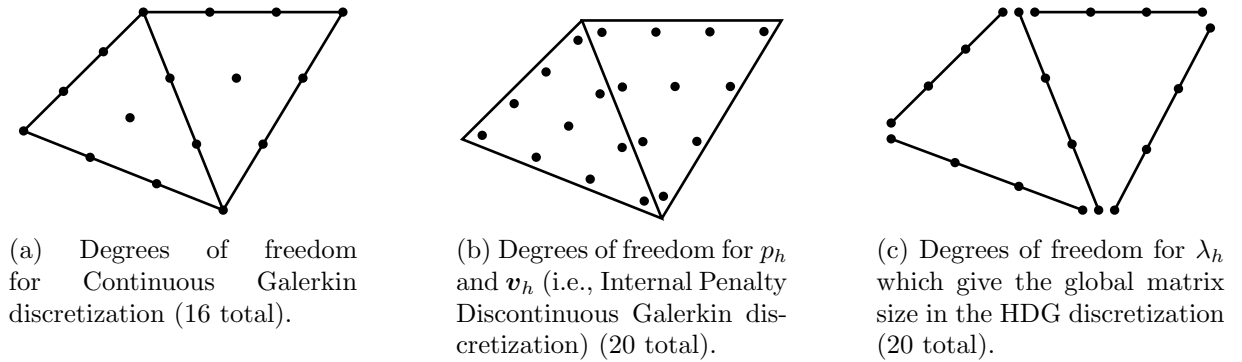

We consider that the approximated pressure and velocity 
are represented in a basis of $\mathbb{P}_{\mathfrak{p}_e}$ 
(e.g., with the Lagrange basis functions) such that, on
every cell, 
\begin{equation}
  \pressure_h^{(e)} = \sum_{k=1}^{N_{\text{dof}}^{(e)}} \mathsf{p}^{(e)}_k \, \phi_k(\bx),
\qquad 
  \velocity_{\bullet~h}^{(e)} = \sum_{k=1}^{N_{\text{dof}}^{(e)}} \mathsf{v}_{\bullet,\, k}^{(e)} 
                              \, \phi_k(\bx),
\end{equation}
The number of degrees of freedom, $N_{\text{dof}}^{(e)}$, 
is given from~\cref{eq:polynomial_space} and amounts at a
given order $\mathfrak{p}_e$ to 
$N_{\text{dof}}^{(e)} = (\mathfrak{p}_e+1)(\mathfrak{p}_e+2)(\mathfrak{p}_e+3) / 6$
for three-dimensional simplexes (tetrahedra).
Thanks to the discontinuous basis functions which make
the solution piecewise polynomial, it is easy to allow 
for different order of approximation depending on the 
cell, i.e., $\mathfrak{p}$-adaptivity.
This is one of the advantage of methods in the Discontinuous 
Galerkin family. 
The positions of the volume degrees of freedom are 
illustrated in \cref{fig:dof:dg}.

We concatenate all coefficients $\mathsf{p}^{(e)}_k$ 
and $\mathsf{v}_{\bullet,\, k}^{(e)}$ so that the 
unknowns associated with the cell $\cell$ are represented 
in a vector by
\begin{equation}
  U_e = \Big( \mathsf{p}_1^{(e)} \,\, \mathsf{p}_2^{(e)} \,\, \ldots \,\,
              \mathsf{p}_{N_\text{dof}^{(e)}}^{(e)} \,\, 
              \mathsf{v}_{x, \, 1}^{(e)} \,\, \ldots \,\,
              \mathsf{v}_{z, \, N_\text{dof}^{(e)}}^{(e)}
\Big)^T.
\end{equation}
Similarly, $\lambda_h$ is represented with a polynomial
on each face, see \cref{fig:dof:hdg}, such that 
\begin{equation}
   \lambda_h \mid_\face ~=~ \sum_{k=1}^{\widehat{N}_{\text{dof}}^{(\face)}} 
                        \lambda_{h,k}^{(\face)} \,\,\, \xi_k(\bx) \, , \qquad
                        \forall \, \face \in \facesetall .
\end{equation}
where $\widehat{N}_{\text{dof}}^{(\face)}$ refers to the number
of degrees of freedom for the face $\face$. 
We define the vector $\Lambda$ to assemble all of the 
coefficients such that
\begin{equation}
  \Lambda = \Big( \, \lambda^{(\face_1)}_{h,1} \, \, \lambda^{(\face_1)}_{h,2} \, \ldots \,
               \lambda^{(\face_\nfaceall)}_{h,\widehat{N}_{\text{dof}}^{(\face_k)}} \, \Big)^T
          = \Big( \, \lambda_{h,1} \, \, \lambda_{h,2} \, \ldots \,
               \lambda_{h,\widehat{N}_{\text{dof}}^{\facesetall}} \, \Big)^T,
\end{equation}
where $\widehat{N}_{\text{dof}}^{\facesetall}$ indicates the total 
number of degrees of freedom for the faces. 
We also introduce the \emph{connectivity map} $\hdgProject_e$, 
which gives, from the global array $\Lambda$, the degrees of 
freedom on the current cell faces such that 
\begin{equation} \label{eq:hdg_connectivity_map}
  \hdgProject_e \Lambda = \Lambda\mid_\dcell .
\end{equation}

Therefore, we have a volume discretization 
for $\pressure_h$ and $\velocity_h$ with
piecewise polynomials on each cell, while $\Lambda$
is defined on the skeleton (piecewise-polynomial on 
the faces) of the complete mesh. 
We illustrate in \cref{fig:dof} with the positions 
of the degrees of freedom. 

We insert the discretized representation in the 
local problem (the two first equations 
of \cref{eq:variational_discretization_full}), 
and have, on each cell,
\begin{equation} \label{eq:hdg_local}
 \hdgA_e U_e + \hdgC_e \hdgProject_e \Lambda = \hdgS_e \, , 
 \quad\qquad \text{HDG discretization: local system.}
\end{equation}

These matrices are defined by testing with 
respect to each function in the basis of 
polynomial. 
The squared matrix $\hdgA_e$ is defined by
\begin{equation} \label{eq:hdg_matrix_A}
   \hdgA_e = \begin{pmatrix}
             -(\freq \kappa^{-1} \phi_i , \, \phi_j)_\cell + \tau (\phi_i,\phi_j)_\dcell
                        & (\partial_x \phi_i , \, \phi_j)_\cell & 
                          (\partial_y \phi_i , \, \phi_j)_\cell & (\partial_z \phi_i , \, \phi_j)_\cell   \\
            -(\phi_i , \, \partial_x \phi_j)_\cell   & -(\freq \rho \phi_i, \phi_j) & 0 & 0       \\
            -(\phi_i , \, \partial_y \phi_j)_\cell   & 0 & -(\freq \rho \phi_i, \phi_j)_\cell & 0 \\
            -(\phi_i , \, \partial_z \phi_j)_\cell   & 0 & 0 & -(\freq \rho \phi_i, \phi_j)_\cell \\
             \end{pmatrix} ,
\end{equation}
where $(\cdot , \, \cdot)$ denotes the inner 
product $( \phi_1 , \, \phi_2)_\cell = 
\int_\cell \phi_1 \, \overline{\phi_2} \, \dd \cell$.
Here, the index $i$ changes with the line
while $j$ changes with the column, both take
values from $1$ to $N_\text{dof}^{(e)}$.
The matrix $\hdgC_e$ 
and the right-hand side $\hdgS_e$ 
are given by
\begin{equation} \label{eq:hdg_matrix_C_S}
 \hdgC_e = \begin{pmatrix}
   - \tau (\xi_k, \, \phi_j)_{\face_1} &  \ldots &    - \tau (\xi_k, \, \phi_j)_{\face_{N_\text{face}}} \\
    (\xi_k, \, \phi_j \n_x)_{\face_1} &  \ldots & (\xi_k, \, \phi_j \n_x)_{\face_{N_\text{face}}}       \\
    (\xi_k, \, \phi_j \n_x)_{\face_1} &  \ldots & (\xi_k, \, \phi_j \n_y)_{\face_{N_\text{face}}}       \\
    (\xi_k, \, \phi_j \n_x)_{\face_1} &  \ldots & (\xi_k, \, \phi_j \n_z)_{\face_{N_\text{face}}}
           \end{pmatrix} \, ; \qquad 
   \hdgS_e = \begin{pmatrix}
                 (f, \, \phi_j)_\cell \\
                    0 \\
                    0 \\
                    0
             \end{pmatrix} \, ,
\end{equation}
where $N_\text{face}$ indicates the number of faces
for the element (i.e., $N_\text{face}=4$ for tetrahedron).
Here, the index $k$ varies with the column and $j$ with the
line. 
In \cref{table:quantities}, we review the notation and 
give the dimension of the matrices for numerical implementation.

\medskip


We proceed similarly with the conditions derived 
on the faces, that is, \cref{eq:lambda_boundary-element,eq:lambda_boundary-domain},
and sum over the cells.
Using the transposed of the connectivity matrix 
$\hdgProject$ of~\cref{eq:hdg_connectivity_map}
to convert the local index of the degrees 
of freedom on the face to their global ones, we have
\begin{equation} \label{eq:hdg_continuity}
 \sum_e \hdgProject_e^T \, 
        \bigg( \, \hdgB_e  \, U_e \, +  \, \hdgL_e  \, \hdgProject_e \Lambda \, \bigg)= 0 \, .
\end{equation}
The matrices $\hdgB_e$ and $\hdgL_e$ are composed 
of blocks for each face of the cell, such that
\begin{equation} \label{eq:hdg_matrix_B_L}
  \hdgB_e = \begin{pmatrix}[1.2]
               \hdgB_e^{(\face_1)} \\ 
               \vdots              \\ 
               \hdgB_e^{(\face_{N_\text{face}})}
            \end{pmatrix} \, ; \qquad
  \hdgL_e = \begin{pmatrix}[1.2]
               \hdgL_e^{(\face_1)} & 0 & \ldots & 0 \\
                 &  & \hspace*{-10mm} \ddots  &              \\
               0 & \ldots & 0 & \hdgL_e^{(\face_{N_\text{face}})}
            \end{pmatrix} \, .
\end{equation}
The faces of each cell can either be an interior 
one or a boundary one, that we respectively denote
by $\face_i^{I}$ and $\face_i^{B}$. 
The matrix $\hdgB_e$ remains the same in all cases:
for the face $k$, we have the matrix 
(with the number of lines given by the number of
 degrees of freedom on the faces, see~\cref{table:quantities})
\begin{equation} \label{eq:hdg_matrix_B_block}
   \hdgB_e^{(\face_i)} = \begin{pmatrix}
   \tau_{\face_i} (\xi_k, \, \phi_j )_{\face_i}  &  
                  (\xi_k, \, \phi_j \n_{x, \face_i})_{\face_i}  &
                  (\xi_k, \, \phi_j \n_{y, \face_i})_{\face_i}  &  
                  (\xi_k, \, \phi_j \n_{z, \face_i})_{\face_i}  \\
        \end{pmatrix} \, .
\end{equation}
For the matrix $\hdgL_e$, each block is squared
with, however, a different definition for interior 
and boundary faces, such that
\begin{equation} \label{eq:hdg_matrix_L_block}
  \hdgL_e^{(\face_i^I)} = - \tau_{\face_i^I} (\xi_k, \, \xi_j )_{\face_i^I}
  \, , \qquad \qquad
  \hdgL_e^{(\face_i^B)} = \Big( \dfrac{\alpha}{\freq \rho \beta} - 
                                \tau_{\face_i^B} \Big) \,\,
                          (\xi_k, \, \xi_j )_{\face_i^B} \, .
\end{equation}
We review the quantities in \cref{table:quantities}.

\bigskip

To assemble the global system, we first replace 
$U_e$ in \cref{eq:hdg_continuity} 
using \cref{eq:hdg_local}, we get,
\begin{equation}
 \sum_e \hdgProject_e^T \, 
        \bigg(
  \, \hdgB_e  \hdgA_e^{-1} \Big( \hdgS_e - \hdgC_e \hdgProject_e \Lambda \Big)
  \, +  \, \hdgL_e  \, \hdgProject_e \Lambda \, \bigg) = \, 0 \, .
\end{equation}
After rearrangement, we obtain the global system 
for the HDG discretization:
\begin{equation} \label{eq:hdg_global}
 \sum_e \hdgProject_e^T 
        \Big( \hdgL_e \, -  \, \hdgB_e  \, \hdgA_e^{-1} \hdgC_e \,  \Big) 
        \, \hdgProject_e \Lambda \, = \, 
      - \sum_e \, \hdgProject_e^T \, \hdgB_e \, \hdgA_e^{-1} \, \hdgS_e \, .
\end{equation}

\begin{table}[ht!] \begin{center}
\caption{Summary of numerical quantities to implement the HDG 
         discretization. The (sparse) global linear system is 
         of size $\widehat{N}_{\text{dof}}^{\Sigma}$ to 
         retrieve the coefficients of $\Lambda$; it is followed 
         by (dense) local linear systems of smaller size on each cell.}
\label{table:quantities}
\renewcommand{\arraystretch}{1.45}
\begin{tabular}{|>{\centering\arraybackslash}p{.07\linewidth}|
                 >{\arraybackslash}p{.82\linewidth}|}
\hline
$\mathsf{dim}$          & Problem dimension.                        \\ \hline
$N_\text{face}$         & Number of face of an element ($N_\text{face}=4$ for 3D tetrahedron). \\ \hline
$N_{\text{dof}}^{(e)}$  & Number of volume degrees of freedom on a cell $\cell$ (for $\pressure_h$ and $\velocity_{\bullet,h}$).\\ \hline
$\widehat{N}_{\text{dof}}^{(\face)}$ & Number of degrees of freedom on the face $\face$ (for $\lambda_h$). \\ \hline
$\widehat{N}_{\text{dof}}^{\Sigma}$  & 
                        Total number of face degrees of freedom, 
                        $\widehat{N}_{\text{dof}}^{\Sigma} = \sum_\face \widehat{N}_{\text{dof}}^{(\face)}$. \\ \hline
$\Lambda$               & Vector of size $\widehat{N}_{\text{dof}}^{\Sigma}$ 
                          encompassing all coefficients of the $\lambda_h$. \\ \hline
$U_e$                   & Vector of size $N_{\text{dof}}^{(e)}$ 
                          encompassing the coefficients of $\pressure_h$
                          and $\velocity_{\bullet,h}$ on the cell $\cell$. \\ \hline
$\hdgA_e$               & Matrix of size 
  $\big( (\mathsf{dim}+1) \, N_{\text{dof}}^{(e)}\big) \, \times \, \big( (\mathsf{dim}+1) \, N_{\text{dof}}^{(e)} \big)$\,,
                        see \cref{eq:hdg_matrix_A}.   \\ \hline
$\hdgC_e$               & Matrix of size
  $\big( (\mathsf{dim}+1) \, N_{\text{dof}}^{(e)}\big) \, \times \, 
   \big( \sum_{\face\in\dcell} \widehat{N}_{\text{dof}}^{(\face)}) \big)$\,, 
                        see \cref{eq:hdg_matrix_C_S}. \\ \hline
$\hdgS_e$               & Right-hand side vector of size $(\mathsf{dim}+1) \, N_{\text{dof}}^{(e)}$\,, 
                        see \cref{eq:hdg_matrix_C_S}. \\ \hline
$\hdgB_e$               & Matrix of size
  $\big( \sum_{\face\in\dcell} \widehat{N}_{\text{dof}}^{(\face)}) \big) \, \times \, 
   \big( (\mathsf{dim}+1) \, N_{\text{dof}}^{(e)}\big)$\,,
                        see \cref{eq:hdg_matrix_B_L,eq:hdg_matrix_B_block}. \\ \hline
$\hdgL_e$               & Matrix of size
  $\big( \sum_{\face\in\dcell} \widehat{N}_{\text{dof}}^{(\face)}) \big) \, \times \, 
   \big( \sum_{\face\in\dcell} \widehat{N}_{\text{dof}}^{(\face)}) \big)$\,,
                        see \cref{eq:hdg_matrix_B_L,eq:hdg_matrix_L_block}. \\ \hline
\end{tabular} \end{center}
\end{table}

\begin{remark}[Boundary conditions] \label{rk:hdg-boundary-conditions}
  If one considers Neumann boundary conditions for \cref{eq:euler_main_BC},
  it amounts to taking $\alpha = 0$ such that  there is no more distinction
  between interior and exterior boundaries in \cref{eq:hdg_matrix_L_block}. 
  However, for Dirichlet boundary condition ($\beta = 0$), one cannot use 
  formula \cref{eq:hdg_matrix_L_block} due to the singularity.
  In this case, the solution is to consider an identity block for $\hdgL_e^{(\face_i^B)}$, 
  which in turn enforces that the values of the trace is zero,
  while we also impose $\hdgB_e^{(\face_i^B)}=0$. \hfill $\triangle$
\end{remark}

\begin{remark}[Numerical approximation of the integrals] \label{remark:integral}
  In the numerical implementations, there are commonly 
  two ways to approximate the value of the integrals 
  that are needed to form the matrices. 
  On the one hand, one can use the quadrature rules 
  for the integration of polynomial functions. 
  On the other hand, one can work with a \emph{reference element},
  that is, find the geometrical transformation from an
  arbitrary simplex to a regular one, and then
  explicitly obtain the integral of the polynomials.
  Here, we use quadrature rules, which are  
  stable for high-order polynomials, and allow to 
  easily consider model parameters that vary within a cell.
\end{remark}

\subsection{Numerical features}
\label{subsection:numerical_features}

The numerical discretization using HDG follows two 
layers, with a global and local linear systems to 
be solved, respectively~\cref{eq:hdg_global,eq:hdg_local}.
The global linear system is written with respect to
the degrees of freedom (dof) of $\lambda_h$, such
that its size is the total number of face dof only
(dimension of $\Lambda$).
Therefore, contrary to CG and DG, HDG avoids the inner cell 
dof for the global problem, hence reducing the size of 
the linear system to be solved, upon taking a sufficiently
high order. 
This is particularly useful for applications such as 
seismic or helioseismology, where the time-harmonic 
approach remains overwhelming for large domains.

Next, the local systems are a specificity of HDG. 
While it might appear as an overhead burden compared 
to other discretization scheme, it is important to 
note that the local problems~\cref{eq:hdg_local} are 
\emph{independent} by cell.
That means it is `embarrassingly parallelizable', 
i.e., it does not need any communication between 
the working processors. 
In addition, these local systems are usually 
composed of small matrices, as their size is the 
number of dof in the current cells; we illustrate
in \cref{section:perspectives}.
Note also that, similarly to the other methods in the
DG-family, we have independent contributions of each 
cell to the global matrix, which makes it convenient 
in a parallel implementation.

The HDG discretization works with 
the first-order formulation such that
both the scalar pressure field and the vectorial 
velocity are computed \emph{with the same accuracy}
while the global matrix is \emph{only} assembled
for \emph{one, scalar} unknown ($\Lambda$). 
This is another advantage of HDG: in the case of 
other discretizations (such as CG, FD or DG), 
the global matrix works directly with the fields of 
interest such that, if one wants to discretize the first-order
formulation, the size of the system contains both the 
scalar \emph{and} vectorial unknowns
(i.e., in three-dimensions, it means four set of global 
unknowns instead of one for HDG).
As an alternative, one can solve for the scalar unknown only
(the pressure field), using the second-order 
formulation~\cref{eq:euler_pressure} and then deduce the 
vectorial velocity. However this means that the numerical 
approximation for the velocity looses one order of 
accuracy compared to the pressure, 
because of the derivative in~\cref{eq:euler_main}. 
Therefore, one would need more advanced techniques to obtain
both fields with similar accuracy, while this is natural with
the HDG discretization.
It has motivated its used in applications 
where both the velocity and the pressure fields are employed,
cf. \cite{Faucher2019FRgWI,Faucher2020DAS}.

\begin{remark}[$\mathfrak{p}$-adaptivity]
  \label{remark:p-adaptivity}
  Similarly to other discretization methods in the DG family, 
  HDG can easily account for different polynomial orders among 
  the cells. 
  Indeed, one simply has to carefully compute the size of the 
  local matrices. 
  \update{Regarding the order of the polynomial for the traces 
          $\lambda_h$, it is convenient 
          to take the polynomial order on each face as the maximum
          order between the corresponding two adjacent cells.
          This is our choice for the numerical experiments below.} \hfill $\triangle$
\end{remark}

\section{Adjoint-state method for HDG discretization}
\label{section:adjoint_state_hdg}

To perform the iterative minimization for the quantitative 
reconstruction of the model parameters, the derivative of 
the misfit functional \cref{eq:misfit} must be evaluated.
For large-scale applications, the adjoint-state method 
(\cite{Chavent1974}) is the natural choice as it avoids 
the explicit formation of the Fr\'echet derivative $D\forward$.
The method is well-known, e.g., \cite{Plessix2006,Chavent2010,Pham2018,Faucher2020Geo}
but requires some careful steps in the context of HDG, because
we have two levels of discretization (local and global),
contrary to the other discretization approaches.

In the following, for simplicity, we assume that there 
is only one source and one frequency for the misfit 
functional~\cref{eq:misfit}. By linearity, they can be
reintroduced later on. Furthermore, we write in the discretized 
settings,
\begin{equation}
  \misfit(\model) \, = \, \dfrac{1}{2} \,
  \big\Vert \, \mathfrak{R} \, U(\model)  \, - \, \data \, \big\Vert^2_2 \, ,
\end{equation}
where $\mathfrak{R}$ is a restriction operator 
(linear) that maps the numerical solution to
the values at the receivers location.


We follow the steps of the adjoint-state method, and 
write the minimization problem with constraints, 
\begin{equation} \label{eq:misfit-constraints}
  \min_{\model} \, \misfit(\model) \, ,\qquad \text{subject to \cref{eq:hdg_local,eq:hdg_continuity}}.
\end{equation}
The first step is to write the formulation with Lagrangians 
(\cite{Lions1971,Glowinski1985}) and to explicit the constraints:
\begin{empheq}[]{align} \label{eq:main-lagrangian}\begin{split}
  \lagrangian(\model,\, \tildu,\, \tildL,\, \tildadjA, \, \tildadjB) = 
      \misfit(\model) + \sum_e \,\, &
      \big< \hdgA_e \tildu_e + \hdgC_e \hdgProject_e \tildL - \hdgS_e, \,
            \tildadjA \big> \\
 &+   \big< \hdgB_e \tildu_e + \hdgL_e \hdgProject_e \tildL, \,
            \hdgProject_e \, \tildadjB \big>,
\end{split} \end{empheq}
\update{where the last two terms correspond to the weak 
formulation of~\cref{eq:hdg_local,eq:hdg_continuity}.
Here, $\langle \cdot , \, \cdot \rangle$ denotes the complex 
inner product in $L^2$ such that 
$\langle u , \, v \rangle = u^* \, v$, and $^*$ is the adjoint.
The formulation contains two Lagrange multipliers: 
$\tildadjA$ has the same dimension as $U_e$ and 
$\tildadjB$ has the same dimension as $\Lambda$.
Moreover, as \cref{eq:main-lagrangian} is linear 
with respect to $\tildadjA$ and $\tildadjB$, it
is a saddle point problem.}

The derivative of $\lagrangian$ with respect to the model parameter 
$\model$ is
\begin{equation} \label{eq:lagrange_gradient_00}
  \partial_{\model} \bigg(\lagrangian(\model,\, \tildu,\, \tildL,\, \tildadjA, \, \tildadjB) \bigg) \, = \,
  \Real \, \bigg( 
  \dfrac{d \, \lagrangian}{d\, \model} \, + \, 
  \dfrac{\partial \lagrangian}{\partial \tildu} \dfrac{\partial \tildu}{\partial \model} \, + \,   
  \dfrac{\partial \lagrangian}{\partial \tildL} \dfrac{\partial \tildL}{\partial \model}  \, \bigg),
\end{equation}
where we follow \cite[Appendix~A]{Faucher2019IP} for the 
specificity of the derivative with complex-variables.
\update{Upon taking $U$ and $\Lambda$ solutions 
to \cref{eq:hdg_local,eq:hdg_continuity}, that is,
they fulfill the constraints in \cref{eq:misfit-constraints},
we have
\begin{equation} \label{eq:lagrange_gradient_0}
  \partial_{\model} \lagrangian(\model, \, U, \, \Lambda, \, \tildadjA, \, \tildadjB)
            = \nabla_{\model} \misfit(\model) \, . 
\end{equation}

\medskip

The \emph{adjoint states} $\adjA$ and $\adjB$ are 
selected such that the derivatives of the Lagrangian 
$\lagrangian$ with respect to $\tildu$ and $\tildL$ 
equate zero.
Therefore, combining 
\cref{eq:lagrange_gradient_00,eq:lagrange_gradient_0}
with these adjoint variables, we have
\begin{equation} \label{eq:lagrange_gradient_0a}
  \partial_{\model} \lagrangian(\model, \, U, \, \Lambda, \, \adjA, \, \adjB)
            \, = \, \nabla_{\model} \misfit(\model)
            \, = \, \Real \, \bigg( 
                             \dfrac{d \, \lagrangian}{d\, \model} \bigg) .
\end{equation}

The derivative of the Lagrangian~\cref{eq:main-lagrangian} 
with respect to $\tildu_e$ is given by, for all cell $\cell$,
\begin{equation}
    \partial_{\tildu_e} \lagrangian \, = \, \partial_{\tildu_e} \misfit
                                \, + \, 
                                \hdgA^*_e \, \tildadjA 
                                \, + \,
                                \hdgB^*_e \, \hdgProject_e \tildadjB
                                \, ,
\end{equation}
and the derivative with respect to $\tildL$ is
\begin{equation}
  \partial_{\tildL} \lagrangian \, = \,
  \sum_e \Big( \, \hdgProject^T_e \, \hdgC_e^* \, \tildadjA 
               \, + \, \hdgProject^T_e \, \hdgL_e^* \, \hdgProject_e \tildadjB \, \Big) \, ,
\end{equation}
where, because $\hdgProject_e$ is real, 
$\hdgProject_e^* = \hdgProject_e^T$, with $^T$ 
the transposed ($\hdgProject_e$ only 
converts the local indexes to the global ones, 
cf. \cref{eq:hdg_connectivity_map}).
By defining $\adjA$ and $\adjB$ such that the derivatives 
of $\lagrangian$ with respect to $\tildu$ and $\tildL$ 
are zero, they solve the system
\begin{subequations}  \label{eq:adjoint_state_0}
\begin{empheq}[left={\empheqlbrace}]{align}
  \sum_e  \, \, \hdgA_e^* \, \adjA + \hdgB_e^* \, \hdgProject_e \adjB & 
             \, = \, - \, \mathfrak{R}^* \big( \mathfrak{R} \, \tildu  \, - \, \data \, \big) \, , 
                          \label{eq:adjoint_state_0a} \\
  \sum_e \, \, \hdgProject^T_e \, \hdgC_e^* \, \adjA 
             \, + \, \hdgProject^T_e \, \hdgL_e^* \, \hdgProject_e \adjB &= 0 \, .\label{eq:adjoint_state_0b}
\end{empheq} \end{subequations}
This can be written in matrix form such that,
for each cell $\cell$ we have,
\begin{equation}
   \begin{pmatrix}[1.25]
      \hdgA_e^*                     & \qquad \hdgB_e^* \\
       \hdgProject^T_e \, \hdgC_e^* & \qquad \hdgProject^T_e \, \hdgL_e^*
   \end{pmatrix} \begin{pmatrix}[1.25]
   \adjA \\
   \hdgProject_e \adjB
   \end{pmatrix} \, = \, \begin{pmatrix}[1.25]
   \big[ \mathfrak{R}^* ( \mathfrak{R} \, \tildu  \, - \, \data \, ) \big]_e \\
   0
   \end{pmatrix} \, , 
\end{equation}
where the notation $[\mathfrak{R}^* ( \mathfrak{R}  \tildu - \data)]_e$
indicates the restriction to the degrees of freedom of the cell $\cell$. 

Similar to the forward problem, we obtain the 
\emph{global} system associated with the adjoint
states by replacing $\adjA$ in~\cref{eq:adjoint_state_0b} 
by its expression from~\cref{eq:adjoint_state_0a}, such that 
\begin{equation}
    \sum_e   \hdgProject^T_e  \hdgC_e^* 
                 \big(-\hdgA_e^{-*} \hdgB_e^*  \hdgProject_e \adjB \big)    
              +  \hdgProject^T_e  \hdgL_e^*  \hdgProject_e \adjB 
              =  \sum_e 
             \hdgProject^T_e  \hdgC_e^*  \hdgA_e^{-*} 
             \big[  \mathfrak{R}^* ( \mathfrak{R}  \tildu   -  \data  )  \big]_e , 
\end{equation}
This can be rearranged as}
\begin{equation}\label{eq:hdg_adjoint_global}
    \sum_e  \hdgProject^T_e  \bigg(  
           \hdgL_e^* - \hdgC_e^*  \hdgA_e^{-*} \hdgB_e^* \bigg) \hdgProject_e \adjB 
           = \sum_e 
             \hdgProject^T_e  \hdgC_e^*  \hdgA_e^{-*} 
              \big[  \mathfrak{R}^* ( \mathfrak{R}  \tildu   -  \data  )  \big]_e  \, .
\end{equation}
We recognize on the left-hand side the adjoint of the
global problem for the forward propagation~\cref{eq:hdg_global}.
The local problems verify~\cref{eq:adjoint_state_0a}, that we 
recall for convenience:
\begin{equation} \label{eq:hdg_adjoint_local}
 \hdgA_e^* \, \adjA  \, = \, - \hdgB_e^* \, \hdgProject_e \adjB  \, - 
           \, \big[ \, \mathfrak{R}^* ( \mathfrak{R} \, \tildu  \, - \, \data \, ) \, \big]_e\, 
           , \qquad \forall \cell \in \mesh.
\end{equation}
\update{The computational steps for the forward and adjoint-state
problems in the context of HDG discretization are summarized in 
\cref{algorithm:adjoint-forward}.}

\setlength{\stepalgo}{0.25em}
\begin{algorithm}[ht!]
\caption{\update{Computational steps for the forward 
         and adjoint-state problems. The local systems
         are solved for all cells $\cell$ of the mesh.}}
\label{algorithm:adjoint-forward}
\begin{mdframed}[rightmargin=22]
  \makebox[0.40\textwidth][l]{compute the global matrix}
        $~\mathcal{A} = \sum_e \hdgProject_e^T \Big( \hdgL_e \, -  
                        \, \hdgB_e  \, \hdgA_e^{-1} \hdgC_e \,  \Big) \, \hdgProject_e$ \\[\stepalgo]
  \textbf{forward problem}  \\
  \hspace*{0.1em} \makebox[0.40\textwidth][l]{compute the forward rhs}
        $\mathcal{B} =  - \sum_e \, \hdgProject_e^T \, \hdgB_e \, \hdgA_e^{-1} \, \hdgS_e$ \\[\stepalgo]
  \hspace*{0.1em} \makebox[0.40\textwidth][l]{solve the gobal system}
        $\mathcal{A} \, \, \Lambda =  \mathcal{B}$ \\[\stepalgo]
  \hspace*{0.1em} \makebox[0.40\textwidth][l]{solve the local systems}
        $U_e \, = \, \hdgA_e^{-1} \Big( - \hdgC_e \hdgProject_e \Lambda + \hdgS_e \Big)$ \\[2\stepalgo]
  \textbf{adjoint-state problem}  \\[1.20\stepalgo]
  \hspace*{0.1em} \makebox[0.40\textwidth][l]{compute the adjoint rhs}
        $\mathcal{C} = \sum_e \hdgProject^T_e  \hdgC_e^*  \hdgA_e^{-*} 
                       \big[  \mathfrak{R}^* ( \mathfrak{R}  U   -  \data  )  \big]_e$ \\[\stepalgo]
  \hspace*{0.1em} \makebox[0.40\textwidth][l]{solve the gobal system}
        $\mathcal{A}^* \, \, \adjB =  \mathcal{C}$ \\[\stepalgo]
  \hspace*{0.1em} \makebox[0.40\textwidth][l]{solve the local systems}
        $\adjA  = - \hdgA_e^{-*} \big(\hdgB_e^* \hdgProject_e \adjB +
                  [\mathfrak{R}^* ( \mathfrak{R} U -\data )]_e\big)$
  \\[\stepalgo]
\end{mdframed}
\end{algorithm}

\medskip

\update{Once the forward and adjoint problems 
have been solved, the gradient of the misfit functional 
is obtained from~\cref{eq:main-lagrangian,eq:lagrange_gradient_0a},}
such that,
\begin{empheq}[]{align} \begin{split}
  \nabla_{\model} \misfit  \, = \,  \Real \bigg( \sum_e
    & \big< (\partial_{\model} \hdgA_e) U_e + (\partial_{\model} \hdgC_e) \hdgProject_e \Lambda, \,
            \adjA \big> \\
 & \qquad + \, \big< (\partial_{\model} \hdgB_e) U_e + (\partial_{\model} \hdgL_e) \hdgProject_e \Lambda, \,
            \hdgProject_e \, \adjB \big> \bigg).
\end{split} \end{empheq}
In fact, for the acoustic equations, only $\hdgA_e$ depends 
on the medium parameter if we ignore the Robin boundary 
conditions which involves $\rho$
(note that in elasticity using the same convention, the matrix
$\hdgC_e$ also depends on the medium parameters, cf. \cite{Bonnasse2017}).
Then, the derivative of $\hdgA_e$ is straightforward
from its definition in~\cref{eq:hdg_matrix_A}.
\medskip

It is crucial that the global problem for the adjoint state
problem is the adjoint of the forward one (with different right-hand side),
see \cref{algorithm:adjoint-forward}, to avoid the re-factorization 
of the matrix. 
The time (and memory) consuming part in 
time-harmonic applications is the matrix 
factorization which, however, allows for the 
fast resolution of multiple rhs problems. 
Here, the factorization of the forward problem can be reused,
using direct solvers  such as \textsc{Mumps} \cite{Amestoy2001,Amestoy2006},
to compute the gradient avoiding an additional matrix factorization
(this is standard with the adjoint-state method for time-harmonic 
equations).
However, we see that the local problems for the forward 
problem are \emph{not} similar to the ones for the adjoint-state
computation. 
Therefore, the code must be adapted with the appropriate 
operations (now using the matrix $\hdgB_e$), but it remains cheap 
thanks to the parallelizability of the local problems,
as discussed in \cref{subsection:numerical_features}.

\medskip


Compared to the adjoint-state method derived for 
the usual discretization (e.g., \cite{Plessix2006,Faucher2019IP}),
the right-hand side is \emph{not} simply made of the 
residuals (the difference between the observations and the simulations), 
but includes the discretization matrices $\hdgA_e$ 
and $\hdgC_e$, see \cref{eq:hdg_adjoint_global}.
Furthermore, note that this modification of the rhs is 
different from the one applied in the forward problem 
(which uses $\hdgB_e$, see~\cref{algorithm:adjoint-forward}).
The local problem for the adjoint involves the residuals, 
and the local matrix $\hdgB_e$ instead of $\hdgC_e$, as 
described in \cref{algorithm:adjoint-forward}.
\section{Three-dimensional experiment}
\label{section:experiments}

\begin{figure}[ht!] \centering
  \subcaptionbox{Three-dimensional domain meshed with about \num{100000} tetrahedron cells.
                \label{fig:fwi:moon-topo_domain_mesh}}
                {\includegraphics[scale=1]{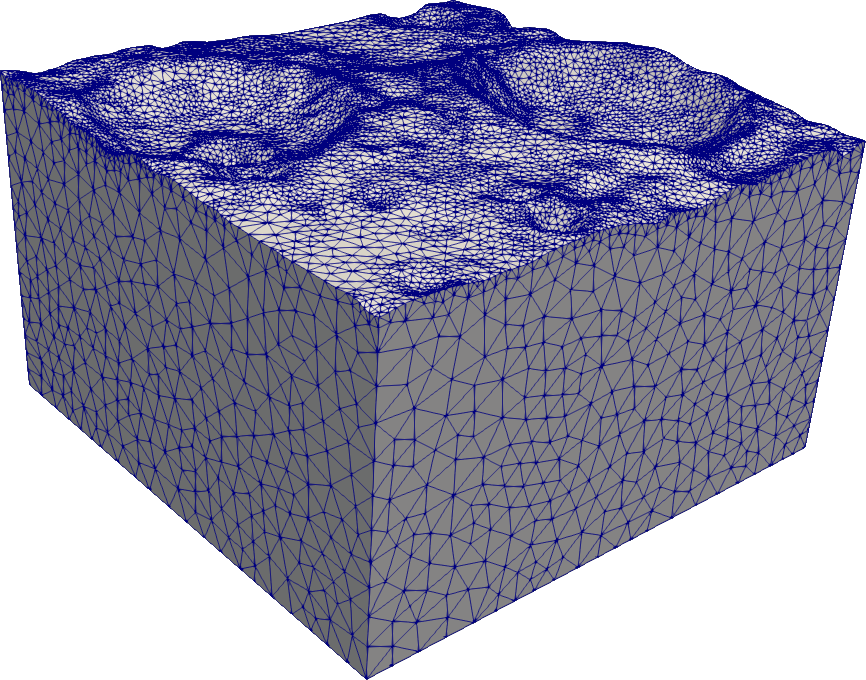}} \hfill
  \subcaptionbox{Cartography of the topography. \label{fig:fwi:moon-topo_domain_topography}}
                {\includegraphics[scale=1]{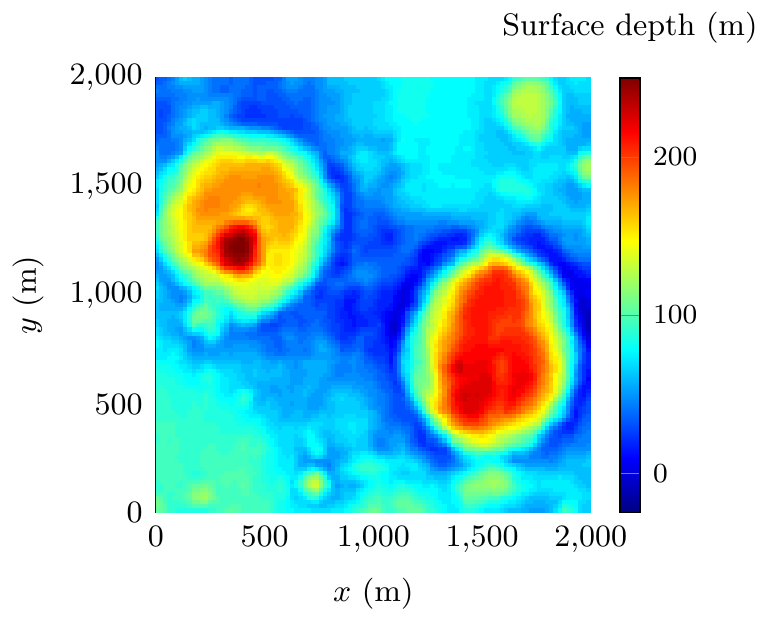}}
  \caption{Domain of interest of size 
           $2 \times 2 \times 1$ \si{\kilo\meter\cubed}, 
           where the surface must be finely meshed to 
           accurately capture the topography. 
           Per convention, the positive values
           indicate the depth.}
  \label{fig:fwi:moon-topo_domain}
\end{figure}

In this section, we illustrate the performance\footnote{
The code we develop and use combines \texttt{mpi} 
and \texttt{OpenMP} for parallelism, it is available 
at \url{https://ffaucher.gitlab.io/hawen-website/}.}
of the HDG discretization in the context of iterative 
reconstruction,
and design a three-dimensional synthetic experiments of size
$2 \times 2 \times 1$ \si{\kilo\meter\cubed}.
The surface is not flat, with a topography made of 
two large craters and smaller variations, it is 
illustrated \cref{fig:fwi:moon-topo_domain_topography}. 
Firstly, it is necessary to accurately capture the topography 
to account for the reflections from the surface, therefore,
we require a very fine mesh of the surface, that we illustrate 
in \cref{fig:fwi:moon-topo_domain_mesh}.
In particular, we rely on the software 
\texttt{mmg}\footnote{\url{https://www.mmgtools.org/}. \label{footnote:mmg}}
to create meshes, and it allows to fix the surface 
cells while possibly coarsening or refining the deeper 
area.
Namely, we use different meshes (to generate the data 
or for the iterations at high-frequency), but we have
the guarantee that the surface remains the same.

\setlength   {\modelwidth}  {6cm}
\setlength   {\modelheight} {6cm}
\renewcommand{\modelfile}   {topo_scale-250_25}
\graphicspath{{figures/fwi/topography/}}

For the numerical computations, the use of 
the HDG discretization is appropriate and allow 
for a flexible framework using 
$\mathfrak{p}$-adaptivity. Indeed, the surface 
is finely mesh, such that low-order polynomials
can be used in the area. On the other hand, in
the deepest part where the cells are larger, we
use higher order polynomials and benefit from the
HDG discretization which disregards the inner dof
for the global linear system.
\update{Following \cref{remark:p-adaptivity}, the 
order for the numerical trace $\lambda_h$ on a given
face is taken to be the maximum order of the two 
adjacent cells that share this face.}

\begin{figure}[ht!] \centering
  \includegraphics[scale=1]{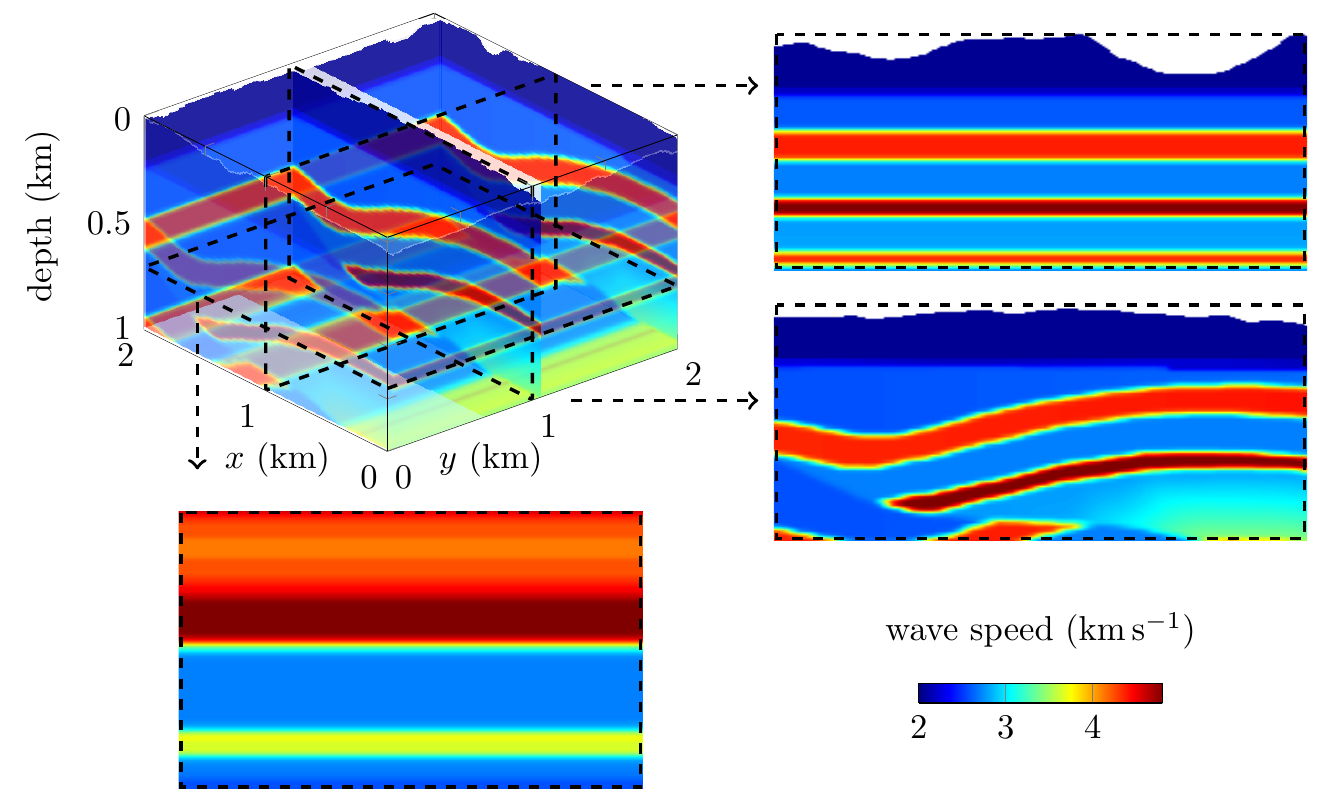}
  \caption{Target wave speed model for the inverse problem, for visualization,
           we extract vertical sections in $x=1$ \si{\km} (top right)
           and $y=1$ \si{\km} (bottom right), and a horizontal section 
           in  $z=\num{500}$ \si{\m} (bottom left).}
  \label{fig:fwi:moon-topo_true}
\end{figure}

In this experiment, the wave speed model contains layers 
of high-contrast velocities, and vary from \num{2000} 
to \num{5500} \si{\m\per\second}, it is show in 
\cref{fig:fwi:moon-topo_true} where we extract 
vertical and horizontal sections for visualization. 
Per simplicity, we consider a constant density with 
$\rho = \num{1000}$ \si{\kg\per\meter\cubed}.
Following a seismic context, a Dirichlet boundary condition
is imposed at the surface (with the topography),
while absorbing boundary conditions are implemented on the other
boundaries.

\setlength   {\modelwidth}  {7cm}
\setlength   {\modelheight} {4cm}
\graphicspath{{figures/fwi/models_scale2000-4800/}}
\renewcommand{\modelfile}   {cp_true}

\subsection{Synthetic partial reflection-data}

We work with reflection data acquired 
from the surface only, with receivers positioned just underneath it
to measure the pressure field. They are positioned to follow the topography
in a lattice with about 75 \si{\meter} between each 
receivers, along the $x$ and $y$ directions. In total, 
we have 625 receivers.
For the data, we consider a set of 100 point-sources 
(i.e., delta-Dirac function for the right-hand side 
of \cref{eq:euler_main}) that are independently excited.
For each of the sources, the 625 receivers measure the 
resulting pressure field. 
The sources are also positioned to follow the topography,
at the surface, in a lattice with about 190 \si{\meter} 
between them, for a total of 100 sources.
This motivates the use of direct solvers which, 
as mentioned, allow for the fast resolution of all sources
once the matrix is factorized. 
As all of the acquisition devices are restricted to the surface
area, the partial data available only consist in \emph{reflection
data}, generated from only one-sided illumination.

We work with synthetic data but include white noise to make 
our experiments more realistic. 
The noise is incorporated in the data with a signal-to-noise 
ratio of $10$ \si{\decibel} for each measurements. 
In addition, the mesh and the order of the polynomial 
for the discretization differ between the generation of 
the data and the inversion procedure.
Following the geophysical setup of this experiment, the 
available frequencies for the reconstruction are limited 
between $5$ and $15$ \si{\Hz}. 
In particular, the absence of low-frequency content in the 
data is an unavoidable difficulty of seismic applications
\cite{Bunks1995,Sirgue2004,Faucher2020Geo}.

\subsection{Iterative reconstruction}

\begin{figure}[ht!] \centering
  \includegraphics[scale=1]{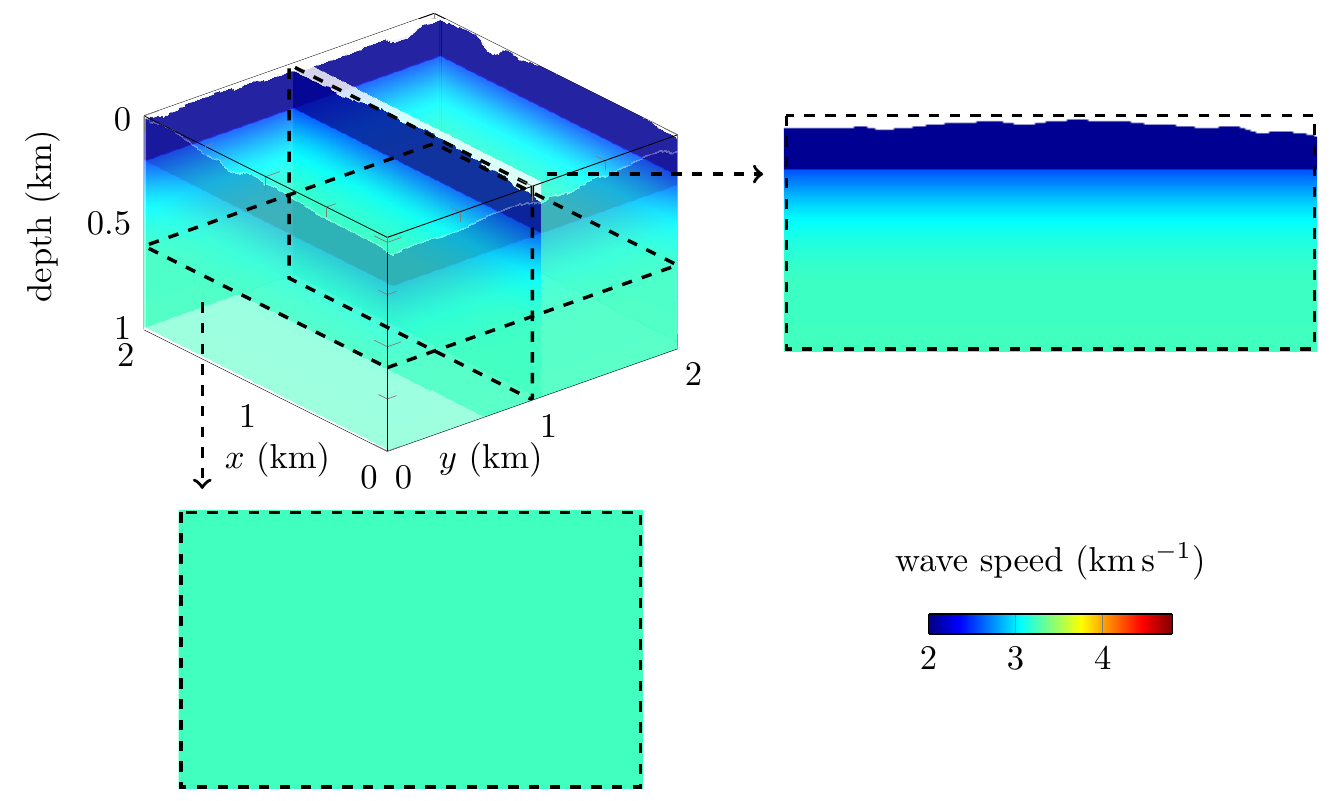}
  \caption{Initial wave speed model for the inverse problem, 
           for visualization, we extract a vertical section
           in $y=1$ \si{\km} (right), and a horizontal section 
           in  $z=\num{500}$ \si{\m} (bottom left).}
  \label{fig:fwi:moon-topo_start}
\end{figure}

We perform the iterative reconstruction using 
data with frequencies from $5$ \si{\Hz} to 
$15$ \si{\Hz}, following a sequential progression,
as advocated in \cite{Faucher2020Geo}. 
The initial model is pictured 
in \cref{fig:fwi:moon-topo_start}: it consists 
in a one-dimensional variation (in depth only)
where none of the sub-surface layers are initially 
known and with an inaccurate background velocity.
In the inverse procedure, the model is represented
as a piecewise-constant for simplicity, that is, 
we have one value of the parameter per cell.

\renewcommand{\modelfile}   {cp_start}

\begin{figure}[ht!] \centering
  \includegraphics[scale=1]{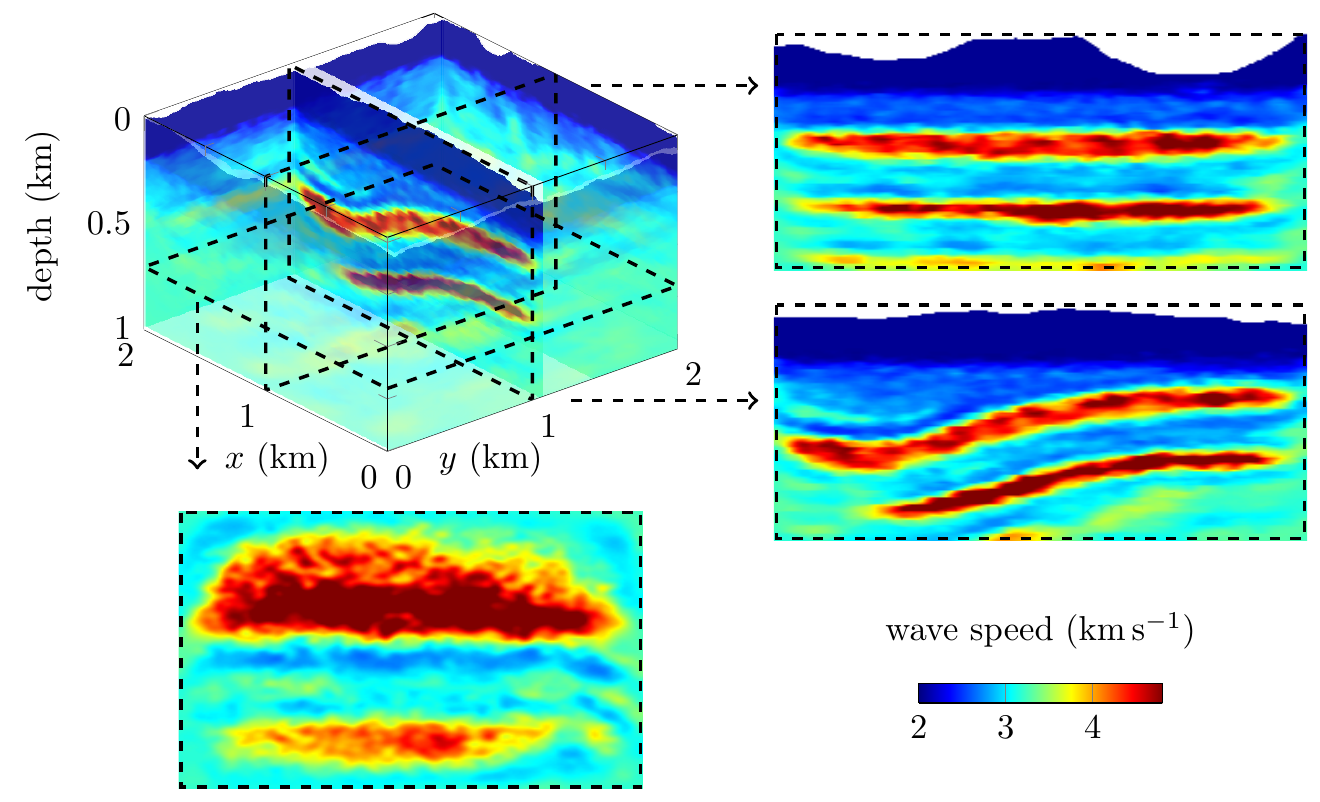}
  \caption{Reconstruction after iterative minimization with frequency 
           up to $15$ \si{\Hz}, for visualization, we extract vertical 
           sections in $x=1$ \si{\km} (top right) and $y=1$ \si{\km} 
           (bottom right), and a horizontal section in $z=\num{500}$ \si{\m} 
           (bottom left).}
  \label{fig:fwi:moon-topo_reconstruction-15hz-p}
\end{figure}

The search direction for the update of the wave speed 
model is computed using the non-linear conjugate 
gradient method, and only depends on the gradient of 
the misfit functional, cf.~\cite{Nocedal2006}.
We perform \num{30} iterations per frequency, for a 
total of \num{300} iterations. 
The order of the polynomials for the basis functions
changes with frequency (increases), while we use two
meshes for the iterations: a coarse mesh for the first
frequencies, from $5$ to $10$ \si{\Hz}, and a more 
refined one for higher frequencies, to capture more
details. 
None of the two meshes for inversion are similar to the one used to
generate the data, but we guarantee that the surface 
cells maintain the same topography.
In \cref{fig:fwi:moon-topo_reconstruction-15hz-p}, we
picture the reconstructed wave speed.

\graphicspath{{figures/fwi/reconstruction/}}
\renewcommand{\modelfile}   {cp_pressure_15hz-480k_gauss2}

We observe that the layers with high velocities 
are appropriately recovered: their positions and 
the values are accurate, except near the 
boundaries, due to the limited illumination. 
On the other hand, the lower values are not retrieved, 
and remain almost similar to the starting model. 
This is most likely due to the lack of background 
information in the starting model, which can only 
be recovered using low-frequency content in the data,
cf.~\cite{Gauthier1986,Luo1991,Faucher2020Geo}.

\section{Perspectives}
\label{section:perspectives}

We have illustrated the use of HDG discretization 
in the context of time-harmonic inverse wave problems. 
The perspective is to handle larger-scale applications, 
using the smaller linear systems provided by the 
method, which only accounts for the degrees of freedom 
(dof) on the faces of the cells. 
Precisely, we have in mind problems made of multiple 
right-hand sides (e.g., with seismic acquisition) where
the bottleneck is usually the memory needed for the 
matrix factorization.
For an efficient implementation of the HDG 
discretization, we have the following requirements.
\begin{enumerate}  \setcounter{enumi}{-1}
  \item One needs to use polynomials of relatively high orders, 
        because the inner-cell dof are omitted in the global 
        linear system, i.e., the mesh must be composed of large
        cells.
  \item On the other hand, large cells cannot be allowed on 
        the whole domain (in order to appropriately represent 
        the geometry) and we rely on $\mathfrak{p}$-adaptivity 
        to adapt the polynomial orders to the size of the cells.
  \item Because we use large cells, the model parameters must 
        be carefully represented.
\end{enumerate}
We provide here a preliminary experiment to illustrate, where 
we multiply by ten the size of the previous model, hence with
a \num{20}$\times$\num{20}$\times$\num{10} \si{\km\cubed} domain.
We generate a mesh of about \num{220 000} cells where, 
similarly as above, the surface is refined to correctly 
represent the topography.

\begin{figure}[ht!] \centering
                {
\begin{tikzpicture}
\node[inner sep=0pt,anchor=south west] (padapt) at (0,0)
     {\includegraphics[scale=0.85]{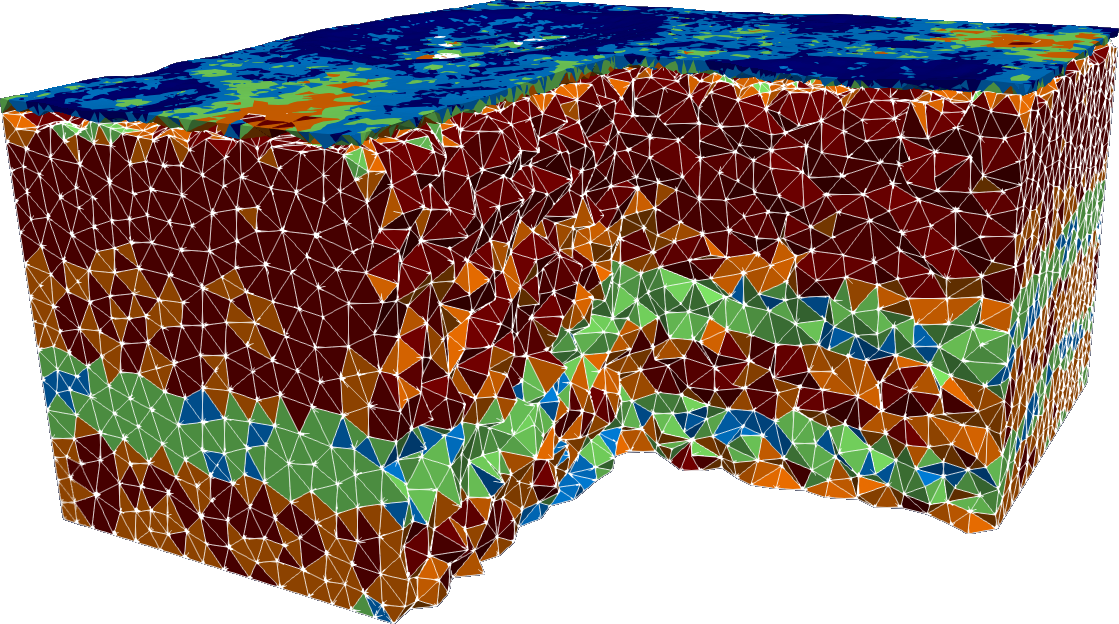}};
      \draw (5,0)  node[anchor=south] {\textbf{(a)}};
\end{tikzpicture}                
      \phantomsubcaption\label{phantom_fig:perspective:P-adaptivity}} \hspace*{-5mm}
{\raisebox{1cm}{
\begin{tikzpicture}
\begin{axis}[height=4.5cm,
             tick label style={font=\footnotesize},
             hide axis,
             colorbar,colormap/jet,
             colorbar style={
             tick label style={font=\scriptsize},
             yshift=0cm,
             width =2mm,
             xshift=-2mm,
             ylabel=\footnotesize{order},
             point meta min=3, point meta max=7,
             }]{}; 
\end{axis}
\end{tikzpicture}
}}\hfill
                {{\raisebox{2.75cm}{
\begin{tikzpicture}
\node[inner sep=0pt,anchor=south west] (padapt) at (0,0)
     {\includegraphics[scale=0.55]{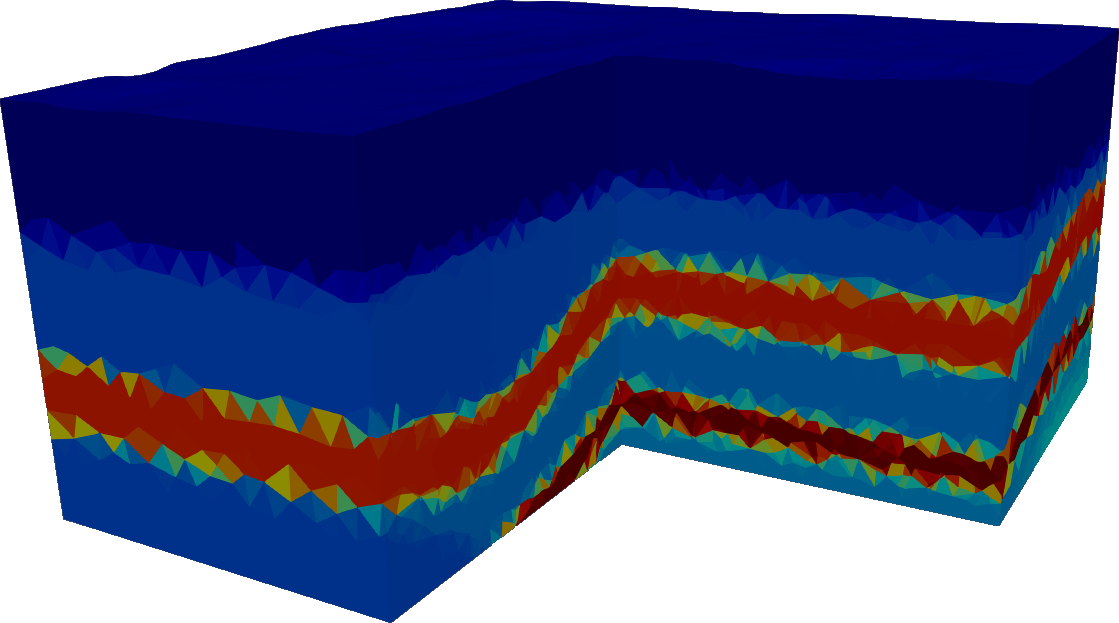}};
      \draw (3.1,0)  node[anchor=south] {\textbf{(b)}};
\end{tikzpicture}
                }}
                 \phantomsubcaption\label{phantom_fig:perspective:P-constant}} \\[-2.20cm]
  \hfill 
                {
\begin{tikzpicture}
\node[inner sep=0pt,anchor=south west] (padapt) at (0,0)
     {\includegraphics[scale=0.55]{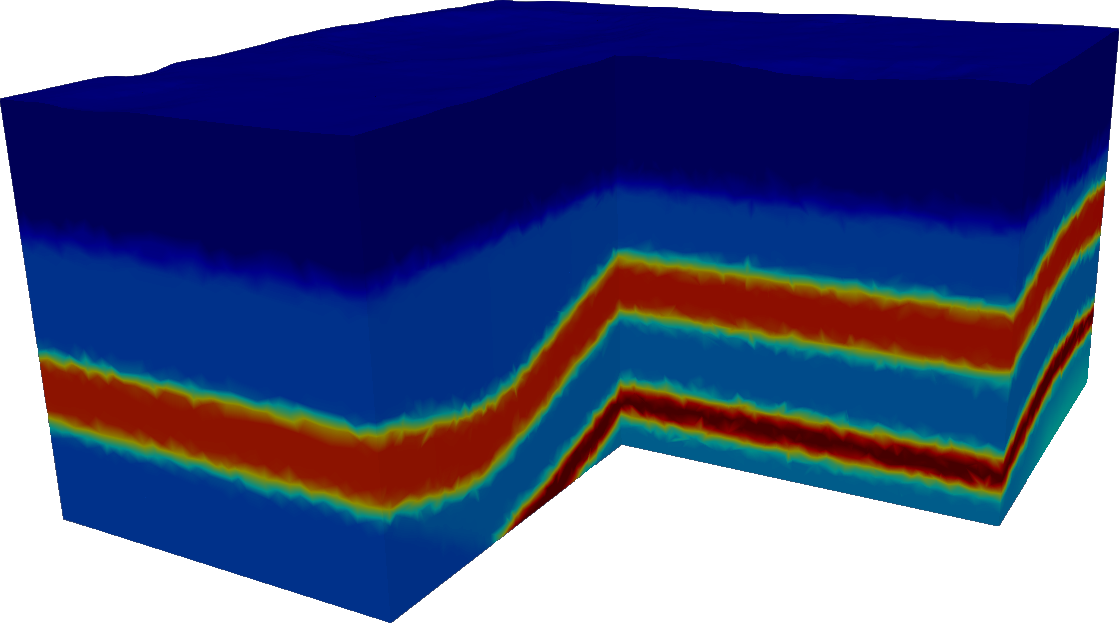}};
      \draw (3.1,0)  node[anchor=south] {\textbf{(c)}};
\end{tikzpicture}
  \phantomsubcaption\label{phantom_fig:perspective:P-poly2}} \\
  \caption{Extension of the model of \cref{fig:fwi:moon-topo_domain}
           to a domain of size 
           \num{20}$\times$\num{20}$\times$\num{10} \si{\km\cubed}.
           \protect\subref{phantom_fig:perspective:P-adaptivity} 
           Illustration of the $\mathfrak{p}$-adaptivity where 
           the order of the polynomial varies per cell depending 
           on the wavelength at $10$ \si{\Hz}.
           \update{The order of the polynomial for the numerical trace 
                   on a face is selected as the maximal order of the 
                   two adjacent cells, see \cref{remark:p-adaptivity}.}
           \protect\subref{phantom_fig:perspective:P-constant} 
           Wave speed model using a piecewise-constant representation
           (one value per cell).
           \protect\subref{phantom_fig:perspective:P-poly2} 
           Wave speed model using piecewise-polynomial representation
           with functions of order 2 per cell.}
  \label{fig:perspective:setup}
\end{figure}

\paragraph{Step 1: $\mathfrak{p}$-adaptivity} 
The first step is to select the order of the polynomials on each cell
(as the dof of each cell are separated, see \ref{fig:dof}).
We rely on the wavelength (that is, the ratio between the frequency and
the wave speed on the cell) to select the order, and we illustrate in
\cref{phantom_fig:perspective:P-adaptivity} the resulting orders at 
$10$~\si{\Hz} frequency.
Here, the order of the polynomials varies from $3$ to $7$, such that the 
cells near surface, which are smaller to handle the topography, consequently 
use a low-order polynomial while the sub-surface layer, made of larger cells 
with low-velocity, need a high order. 
We also observe the sub-surface layer of increasing velocity, 
where the order is allowed to be reduced.\\

\paragraph{Step 2: model representation} 
Because of the large cells, the use of a piecewise-constant model 
parameters with one value per cell can lead to a coarse representation,
that we illustrate in \cref{phantom_fig:perspective:P-constant}. 
Therefore, we rely on a piecewise-polynomial model representation
and illustrate in \cref{phantom_fig:perspective:P-poly2} with 
polynomials of order 2 per cells to represent the wave speed.
Clearly, the latter allows for a precise representation, removing 
the inaccuracies between cells.
In this experiment where the above layer is mostly constant, we can
even envision to use one polynomial on a group of cells instead of 
one per cell. 
Such models can be easily accounted for using the quadrature rules
to evaluate the integrals, see \cref{remark:integral}.

\begin{figure}[ht!] \centering
  \includegraphics[scale=1]{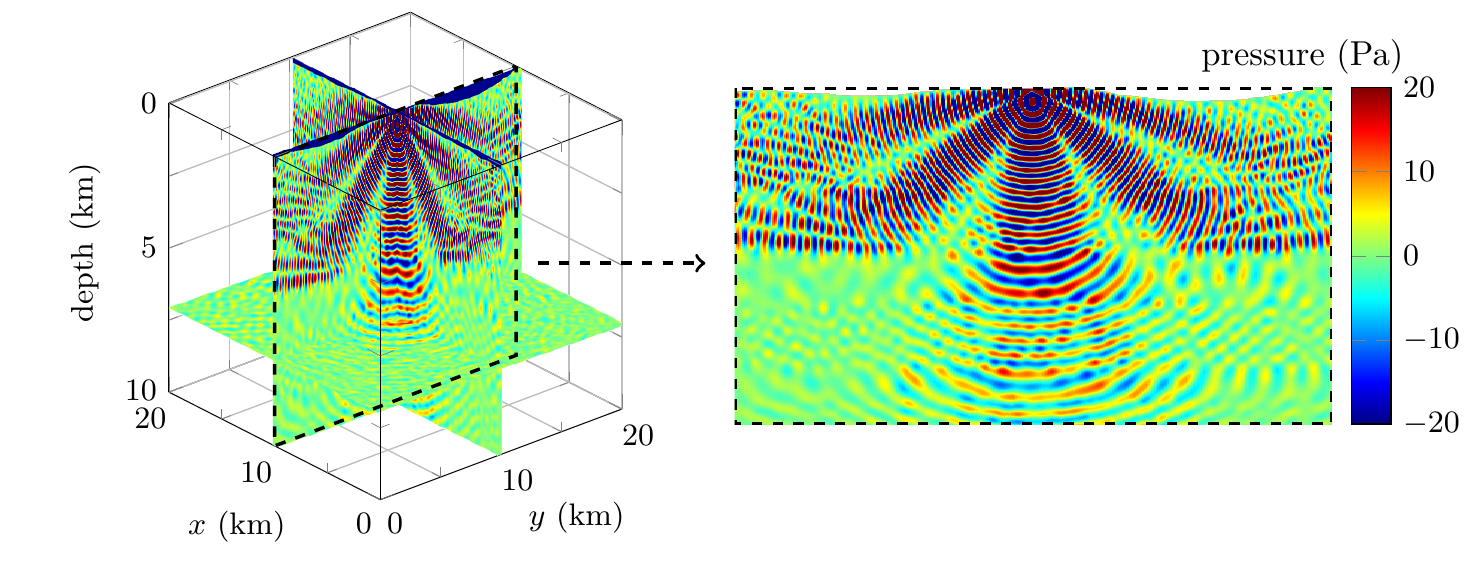}
  \caption{Real part of the time-harmonic pressure 
           field at \num{6} \si{\Hz} on a domain of 
           size \num{20}$\times$\num{20}$\times$\num{10} \si{\km\cubed},
           following the numerical setup given 
           in \cref{fig:perspective:setup}.
           For visualization, we extract a vertical section
           in $x=10$ \si{\km} (right).}
  \label{fig:perspective:wave}
\end{figure}

\paragraph{Numerical cost}
In \cref{fig:perspective:wave}, we illustrate the
solution of the forward problem with the pressure 
field for a single source at $6$ \si{\Hz}, where 
we follow the setup prescribed in \cref{fig:perspective:setup}.
The computation is performed using \num{360} cores 
(using \num{10} nodes with \num{4} processors per node
and \num{9} threads per processor). 
The discretization problem has the following dimensions,
\begin{itemize}
  \item The size of the linear system is 
        $\widehat{N}_{\text{dof}}^{\Sigma} = \num{10993236}$.
  \item The number of volume dof per unknown is 
        $\sum_e N_{\text{dof}}^{(e)} = \num{14812171}$,
        and we have four unknowns, i.e., 
        the pressure and the vectorial velocity.
\end{itemize}
We see that the size of the linear system reduces by \num{25} \%
the number of volume dof for one unknown, and by \num{80} \%
if we consider the total size with the four unknowns.
In this experiment, the matrix factorization requires \num{440} \si{GiB} 
using the direct solver \textsc{Mumps}\footnote{
          Note that the memory cost can be reduced by using 
          different options such as the block low-rank feature 
          recently implemented in the solver, cf. \cite{Amestoy2019}.}.
The time to factorize the global matrix is of 
about \num{6} \si{\min}, while the time to solve 
the local problems to assemble the volume solution 
is less than \num{1} \si{\second}. It highlights 
that the second step of the method where one has 
to solve the local problems is very cheap compared
to the resolution of the global linear system, in 
particular as it does not need any communication 
between the cells.

\begin{remark}
Comparing with other discretization methods, where 
the linear system is made of the volume dof (see \cref{fig:dof}),
we can draw the following remarks,
\begin{itemize} 
  \item Standard DG results in a linear system of size the number of 
        volume dof (see \cref{fig:dof:dg}).
        To solve the first-order system with both the pressure field
        and velocity, it represents an increase of \num{539} \% compared
        to HDG. When solving the second-order system to only recover for
        the pressure field, it is an increase of \num{135} \%.
  \item Using Continuous Galerkin, let us first remind that it is much harder to employ 
        $\mathfrak{p}$-adaptivity because the dof on the faces are shared 
        (see \cref{fig:dof:fe}). 
        Then, assuming a constant order 7, the resulting linear system 
        is increased by \num{490} \% and \num{122} \% compared to HDG, 
        respectively for the first- and second-order systems.
\end{itemize} \hfill $\triangle$
\end{remark}

\section{Conclusion}

We have derived the adjoint-state method for the computation of 
the gradient of a functional in the framework of the Hybridizable 
Discontinuous Galerkin discretization in order to handle large-scale 
time-harmonic inverse problems using direct solvers to account for 
multiple right-hand sides.
HDG reduces the size of the global linear system compared to other 
discretization methods, and it works with the first-order 
formulation of the forward problem.
We have illustrated with the acoustic wave equation, where it gives 
the computational approximation for the pressure and the velocity 
fields at the same accuracy. 
The HDG method relies on two levels, a global system and local ones, 
which must be carefully addressed to obtain the adjoint-state 
where the matrix factorization of the forward problem can still 
be used for the backward problem.
HDG allows to easily account for $\mathfrak{p}$-adaptivity, 
which is useful when some part of the mesh must be particularly 
refined to take into consideration the specificity of the problem, 
such that the geometry of the parameters or, as we have illustrated
in our experiment, with the topography.
We have given a preliminary insight to work with problems of 
larger scales, and we need to continue to investigate if HDG can 
help to fill the gap between the largest time and frequency-domain problems.
Extension to elasticity is straightforward, and requires 
only minor modifications of the steps we have given, large-scale elastic
inversion is part of our ongoing research. 

\section*{Acknowledgments}
 The authors would like to thank
 Ha Pham for thoughtful discussions.
 FF is funded by the Austrian 
 Science Fund (FWF) under the Lise 
 Meitner fellowship M 2791-N.
 OS is supported by the FWF, with SFB F68, project 
 F6807-N36 (Tomography with Uncertainties).

 The code used for the experiments, \texttt{hawen}, is
 developed by FF and is available at 
 \url{https://ffaucher.gitlab.io/hawen-website/}.
 The numerical experiments have been  performed as part 
 of the GENCI resource allocation project AP010411013.

\footnotesize 
\bibliographystyle{siamplain}
\bibliography{bibliography}

\end{document}